\documentclass[preprint,12pt]{elsarticle}
\usepackage{graphicx} 
\usepackage{tabularray}
\usepackage{adjustbox}
\usepackage{amssymb}
\usepackage{amsmath}
\UseTblrLibrary{booktabs}
\usepackage{float}
\usepackage{xcolor}
\usepackage[margin=1.2in]{geometry}
\usepackage{subcaption}
\usepackage{comment}
\usepackage{amsthm}
\theoremstyle{plain}

\newcommand{\ux}{{\mathbf{x}}}

\newcommand{\ua}{{\mathbf{a}}}

\DeclareMathOperator*{\argmin}{argmin}

\makeatletter
\def\ps@pprintTitle{%
   \let\@oddhead\@empty
   \let\@evenhead\@empty
   \def\@oddfoot{\footnotesize\itshape\hfill\today} 
   \let\@evenfoot\@oddfoot}
\makeatother

\begin{document}

\begin{frontmatter}

\title{A Neural Latent Dynamics Approach for Solving Inverse Problems in Cardiac Electrophysiology}

\author[label1]{Edoardo Centofanti\corref{contrib}}
\author[label2]{Giovanni Ziarelli\corref{contrib}}
\author[label3]{Simone Scacchi}
\author[label1]{Luca Franco Pavarino}

\affiliation[label1]{organization={Dipartimento di Matematica, Università di Pavia},
             addressline={Via Adolfo Ferrata, 5},
             city={Pavia},
             postcode={27100},
             country={Italy}
             }
\affiliation[label2]{organization={MOX Laboratory, Dipartimento di Matematica, Politecnico di Milano},
             addressline={Via Edoardo Bonardi, 9},
             city={Milano},
             postcode={20133},
             country={Italy}
             }
\affiliation[label3]{organization={Dipartimento di Matematica, Università di Milano},
             addressline={Via Cesare Saldini, 50},
             city={Milano},
             postcode={20133},
             country={Italy}
             }
\cortext[contrib]{These authors contributed equally}

\begin{abstract}
Solving inverse problems in cardiac electrophysiology consists in the recovery of physiological parameters from surface electrocardiogram (ECG) measurements, a task which is often computationally unfeasible due to the severe ill-posedness and the prohibitive computational complexity of PDE-constrained optimization.
In this work, we introduce a data-driven framework leveraging Latent Dynamics Networks (LDNets) to construct efficient surrogate models of the forward operator.
By mapping low-dimensional parameters, representing ectopic activation sites or ischemic region descriptors, to the ECG signals via latent dynamics governed by neural ordinary differential equations, our approach circumvents the computational burden of evaluating high-fidelity cardiac models during iterative parameter estimation.
The surrogate is trained offline on high-fidelity data, enabling rapid and robust inversion. We validate the proposed framework through rigorous numerical experiments with synthetic data across both 2d and 3d geometries.
Results show that the LDNet-based surrogate achieves precise reconstruction of cardiac parameters while drastically reducing computational overhead, thereby enabling near real-time clinical applications.
\end{abstract}


\end{frontmatter}

\section{Introduction}

The non-invasive identification of cardiac electrical abnormalities from surface electrocardiogram (ECG) signals represents a central challenge in computational cardiology. Clinically relevant applications include the localization of ectopic activation sites responsible for arrhythmias and the detection and characterization of ischemic regions within the myocardium. These tasks can be naturally formulated as inverse problems, where one seeks to infer unknown model parameters or spatial tissue properties from indirect and spatially aggregated measurements \cite{nie07,ruu08,nie09,pul10,ber17,ber20,asp25}.

From a mathematical perspective, the ECG inverse problem is severely ill-posed, due to the smoothing nature of the forward operator, the limited number of measurement locations, and the intrinsic nonlinearity of cardiac electrophysiology models. Small perturbations in the data may lead to large variations in the inferred parameters, and multiple configurations of the cardiac source may produce indistinguishable ECG signals. These issues pose significant challenges for robust and reliable parameter estimation.

A standard approach to tackle this problem relies on high-fidelity (HF) computational models describing the propagation of electrical activity in cardiac tissue. Such models are typically based on the Bidomain or Monodomain equations \cite{col14, sun06}, coupled with detailed ionic models \cite{ten06}, and are solved numerically using advanced discretization techniques such as the finite element method. Within this framework, the inverse problem is commonly formulated as a PDE-constrained optimization problem, in which the forward model is repeatedly evaluated to minimize a discrepancy functional between simulated and observed ECG signals \cite{gra22,gra25}.

While this approach can provide accurate reconstructions, it is computationally extremely demanding. High-fidelity simulations involve the solution of nonlinear, time-dependent PDE systems on fine spatial and temporal grids, often requiring several minutes per simulation even in simplified geometries, and significantly more in three-dimensional settings. As a consequence, their use within iterative optimization procedures, where a large number of forward evaluations is required, becomes prohibitive, especially in multi-query contexts such as uncertainty quantification, parameter calibration, or real-time clinical decision support.

These limitations motivate the development of reduced-order and surrogate models capable of approximating the forward mapping at a drastically reduced computational cost. Classical projection-based techniques, such as Proper Orthogonal Decomposition (POD) \cite{ben15}, have been extensively used to construct low-dimensional representations of parametric dynamical systems. More recently, data-driven approaches based on machine learning have emerged as powerful alternatives, enabling the approximation of complex nonlinear mappings directly from data generated by high-fidelity simulations.

A variety of methodologies have been proposed in this context. Dimensionality reduction techniques, including autoencoders and convolutional autoencoders \cite{wan14}, allow for the construction of compact latent representations of high-dimensional states. The temporal evolution can then be modeled using neural ordinary differential equations \cite{lin22}, recurrent neural networks \cite{mau21}, or other approaches such as Dynamic Mode Decomposition \cite{bru16} and Sparse Identification of Nonlinear Dynamics \cite{sit20}. In parallel, operator learning frameworks, such as Deep Operator Networks \cite{lu21}, Fourier Neural Operators \cite{li20}, and related neural operator architectures \cite{azi24}, have been developed to approximate mappings between infinite-dimensional function spaces.

In this regard, Latent Dynamics Networks (LDNets) \cite{reg24} provide a unified framework that combines dimensionality reduction and dynamical system learning in a single architecture. By learning a low-dimensional latent representation together with its governing dynamics, LDNets enable efficient and accurate reduced-order modeling of complex spatio-temporal processes. Related approaches based on neural model learning have also demonstrated promising results for time-dependent systems \cite{reg19}, including recent applications in computational biology \cite{zia25,zap25,zia26,dok26}.

Despite these advances, the use of machine learning-based surrogate models for inverse problems in cardiac electrophysiology remains relatively limited. In particular, the integration of latent dynamical models within optimization-based inverse frameworks for parameter identification in PDE-governed systems is still largely unexplored.

\begin{figure}[t]
    \centering
    \includegraphics[width=\textwidth]{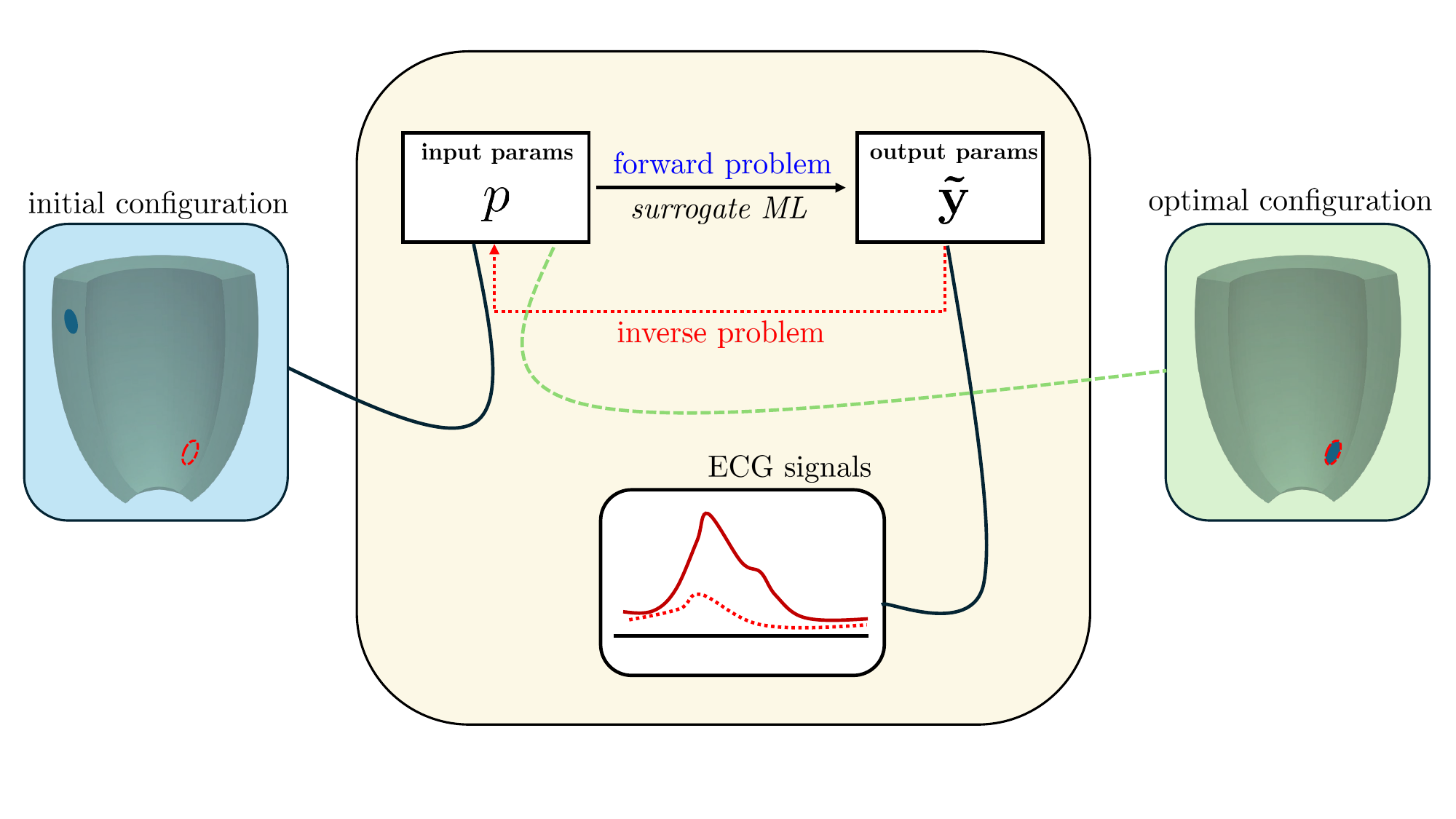}
    \caption{Schematic representation of the approach. Starting from an initial configuration of the input parameters we aim at finding the optimal configuration minimizing the accuracy metric. This measure depends on the discrepancy of ground truth ECG signals and the respective reconstruction. The forward mapping from input parameters to ECG signals is surrogated via the introduced machine-learning strategy.}
    \label{fig:info1}
\end{figure}

In this work, we propose a data-driven framework based on Latent Dynamics Networks to efficiently solve inverse problems arising in cardiac electrophysiology (see Figure \ref{fig:info1} for a schematic representation of the workflow). The key idea is to approximate the forward operator mapping low-dimensional parameters, such as activation site location or ischemic region descriptors, to pseudo-ECG signals, by means of a neural surrogate model defined in a latent space. The surrogate is trained offline using a limited dataset of high-fidelity simulations and subsequently employed within an optimization loop to perform fast parameter estimation.

The proposed approach offers several advantages. First, it significantly reduces the computational cost associated with repeated forward evaluations, enabling efficient exploration of the parameter space. Second, by operating in a low-dimensional parameter setting, it avoids the need to reconstruct the full transmembrane potential field, focusing directly on clinically relevant quantities. Third, the latent dynamics formulation allows for an accurate representation of the temporal structure of ECG signals, which is essential for reliable inversion.

The effectiveness of the proposed framework is demonstrated on several test cases, including the localization of ectopic activation sites and the identification of ischemic regions with fixed and variable extent. Numerical results show that the method achieves accurate reconstructions while maintaining a computational cost compatible with near real-time applications.

The manuscript is organized as follows.
In Section~\ref{sec:mathFrame} we introduce the mathematical framework adopted in this work, detailing the problem formulation, the high-fidelity model, and the surrogate model.
Section~\ref{sec:num_exp} is dedicated to the numerical experiments: we describe the dataset construction and parameter settings in Section~\ref{sec:dataset_details}, present the tailored initialization and optimization strategies in Section~\ref{sec:guess}, and report the results for the different test cases in Sections~\ref{sec:2d}--\ref{sec:2d_isch_rad}.

\section{Mathematical Framework}\label{sec:mathFrame}

In this section, we introduce the core mathematical models underlying the proposed inverse problem framework and describe the associated numerical methods.
We first describe the high-fidelity forward model used to compute pseudo-ECGs, which serves as ground-truth training/testing data for the surrogate model.
Hence, we present the LDNet-based approach and the multiple-shooting  strategy adopted for the inverse problem. The latter helps to accelerate convergence while reducing the risk of getting stuck in undesirable local minima due to the non-convex nature of the problem.

\subsection{Forward Problem: models and methods}

In this section we describe a high-fidelity mathematical model that has been used in literature for reconstructing pseudo-ECG signals at body-surface lead locations on a static torso domain, given an input current defined on the cardiac domain that drives the electrophysiological propagation \cite{sca23}.
In particular, we used this model for generating training, validation and test data.
It consists in solving, in cascade, (a) a PDE model for the electrical potential in a conductive medium and (b) a Laplace problem for recovering the ECG signals. The two steps are described below.

\subsubsection{High-Fidelity Model for the propagation of the electrical signal}

The potential spreading in the cardiac domain can be described by the Monodomain model \cite{ver16}, obtained from the Bidomain equations \cite{afr22, cen24} under the assumption of equal anisotropy ratios between the intra- and extracellular conductivity tensors, namely $\mathbf{D}_e = \lambda \mathbf{D}_i$ with $\lambda \in \mathbb{R}^+$ (cfr., e.g., \cite{col14, sun06}). This assumption leads to a single parabolic partial differential equation coupled with a system of ordinary differential equations governing the local reaction dynamics.

Let $\Omega \subset \mathbb{R}^d$, $d=2,3$, denote the computational cardiac domain and $T>0$ the final simulation time. The Monodomain system reads as:
\begin{equation}
\label{eq:monodomain}
\left\{
\begin{aligned}
&\chi C_m \frac{\partial v}{\partial t}
- \nabla\!\cdot(\mathbf{D}_m \nabla v)
+ I_{\mathrm{ion}}(v,\mathbf{w},\mathbf{c})
= I_{\mathrm{app}}
&& \text{in } \Omega \times (0,T), \\
&\frac{\partial \mathbf{w}}{\partial t}
= \mathbf{R}(v,\mathbf{w})
&& \text{in } \Omega \times (0,T), \\
&\frac{\partial \mathbf{c}}{\partial t}
= \mathbf{C}(v,\mathbf{w},\mathbf{c})
&& \text{in } \Omega \times (0,T), \\
&\mathbf{n}^\top \mathbf{D}_m \nabla v
= 0
&& \text{on } \partial\Omega \times (0,T), \\
&v(\mathbf{x},0)=v_0(\mathbf{x}), \quad
\mathbf{w}(\mathbf{x},0)=\mathbf{w}_0(\mathbf{x}), \quad
\mathbf{c}(\mathbf{x},0)=\mathbf{c}_0(\mathbf{x})
&& \text{in } \Omega.
\end{aligned}
\right.
\end{equation}

Here, $v:\Omega\times(0,T)\to\mathbb{R}$ denotes the transmembrane potential, $\textbf{w}:\Omega\times(0,T)\to\mathbb{R}^{s_w}$ the recovery-gating variables, $\textbf{c}:\Omega\times(0,T)\to\mathbb{R}^{s_w}$ the vector of ionic concentrations, $\mathbf{D}_m$ is the monodomain symmetric positive-definite conductivity tensor, given by 
\[
\textbf{D}_m = \textbf{D}_i (\textbf{D}_i+\textbf{D}_e)^{-1}\textbf{D}_e, 
\]
where $\textbf{D}_{i,e}$ are the intra and extra cellular conductivity tensors. We recall that $\textbf{D}_{i,e}$ are defined as

\begin{equation}
        \textbf{D}_{i,e}(\ux) =                             \sigma_{t}^{i,e} \textbf{I} + (\sigma_l^{i,e}-\sigma_t^{i,e}) \ua_l(\ux)\ua_l^\top(\ux), 
\end{equation}
where $\ua_l$ $\in [L^\infty(\Omega)]^3$ represents the local fiber direction and we have assumed that the tissue is transversely isotropic.
Furthermore, we have denoted as $\sigma_l^{i,e}$ and $\sigma_t^{i,e}$ the conductivity coefficients for 
the intra- and extracellular media along the longitudinal ($\ua_l$) and transversal ($\ua_t$) directions. 

$C_m$ is the membrane capacitance and $\chi$ is the membrane surface to volume ratio.

The nonlinear reaction terms $I_{\mathrm{ion}}$, $\mathbf{R}$, and $\mathbf{C}$ define the local dynamics and depend on the specific choice of ionic or phenomenological model; for the purposes of this work, they are regarded as known nonlinear operators. The applied current $I_{\mathrm{app}}$ acts as a forcing term and represents an external excitation of the system. While the equal anisotropy assumption is a modeling simplification, in the absence of extracellular current injection it provides an accurate and computationally efficient approximation of the Bidomain equations (see, e.g., \cite{col14}).

\subsubsection{High-Fidelity model for computing pseudo-ECGs}

To connect the state variable $v$ to measurable ECG signals, we need to introduce an observation operator.
In our case it will be based on the infinite-volume conductor approximation \cite{sca23}. Indeed, under the assumption that the surrounding medium is isotropic with constant conductivity $\sigma_b$ and extends to infinity, the extra-cardiac potential $u(\mathbf{x}) :\mathbb{R}^3\to\mathbb{R}$ at each time and choice of parameters satisfies
\begin{equation}
\label{eq:extracellular}
-\sigma_b \Delta u(\mathbf{x}) =
\begin{cases}
\nabla\!\cdot(\mathbf{D}_i(\mathbf{x}) \nabla v(\mathbf{x})) & \mathbf{x}\in\Omega, \\
0 & \mathbf{x}\in\mathbb{R}^3\setminus\Omega.
\end{cases}
\end{equation}

Exploiting the fundamental solution of the Laplace operator in three dimensions, the potential recorded at a measurement location $\mathbf{x}\in\mathbb{R}^3\setminus\Omega$, called \textit{pseudo-electrocardiogram} (pseudo-ECG), corresponding to $u(\mathbf{x})$ evaluated at each time for the chosen set of model parameters $p$, can be expressed as 

\begin{equation}
\label{eq:pecg}
\text{pECG}(\mathbf{x},t;p)
\;=\;-
\frac{1}{4\pi\sigma_b}
\int_{\Omega}
\mathbf{D}_i(\mathbf{y}) \nabla v(\mathbf{y},t;p)
\cdot
\nabla\!\left(\frac{1}{|\mathbf{x}-\mathbf{y}|}\right)
\,\mathrm{d}\mathbf{y}.
\end{equation}
Since the locations $\mathbf{x}$ in which we evaluate \eqref{eq:pecg} are fixed, in the following we will refer to $\text{pECG}(\mathbf{x},t;p)$ as $\text{pECG}(t;p): [0,T] \rightarrow\mathbb{R}^{N_{\mathrm{leads}}}$. Further details can be found also in~\cite{col14, ges83}. 

\subsubsection{Numerical solution of the high-fidelity models}\label{sec:num_schemes}

High-fidelity solutions have been computed using the finite element method (FEM), as described in the following.
In particular, we employed Q1 elements, namely square (2d case) or hexahedral (3d case) cells with node at the vertices.
Details on the weak formulation and FEM for the monodomain and the bidomain model, from which the monodomain is derived, can be found in~\cite{col14,cen24}.

For the time evolution, we employed a first order implicit-explicit (IMEX) scheme.
In particular, starting from a standard Galerkin discretization procedure, we can rewrite~\eqref{eq:monodomain} in matrix form as

\begin{equation}\label{eq:mono-matrix}
\begin{cases}
\displaystyle \chi C_m M \frac{d\mathbf{v}}{dt}
+ A \mathbf{v}
+ M I_{\mathrm{ion}}(\mathbf{v},\mathbf{w},\mathbf{c})
= M \mathbf{I}_{\mathrm{app}} \\
\displaystyle \frac{d\mathbf{w}}{dt} = \mathbf{R}(\mathbf{v},\mathbf{w}) \vspace{0.2cm}\\
\displaystyle \frac{d\mathbf{c}}{dt} = \mathbf{C}(\mathbf{v},\mathbf{w},\mathbf{c})
\end{cases}
\end{equation}
where $M$ denotes the finite element mass matrix and $A$ is the stiffness matrix associated with the Monodomain conductivity tensor $\mathbf{D}_m$.
In this case, we define $\mathbf{v}$ as the vector of degrees of freedom associated with the transmembrane
potential $v$, and
\[
\mathbf{R}(\mathbf{v},\mathbf{w})
= (R_1(\mathbf{v},\mathbf{w}),\ldots,R_{s_w}(\mathbf{v},\mathbf{w}))^\top,
\qquad
\mathbf{C}(\mathbf{v},\mathbf{w},\mathbf{c})
= (C_1(\mathbf{v},\mathbf{w},\mathbf{c}),\ldots,C_{s_c}(\mathbf{v},\mathbf{w},\mathbf{c}))^\top.
\]

\noindent
Finally, vectors $\mathbf{I}_{\mathrm{ion}}(\mathbf{v},\mathbf{w},\mathbf{c})$ and
$\mathbf{I}_{\mathrm{app}}$ represent the finite element coefficient vectors
associated with the nonlinear ionic current $I_{\mathrm{ion}}(v,\mathbf{w},\mathbf{c})$
and the applied current $I_{\mathrm{app}}$, respectively.

In the IMEX strategy, we consider the following scheme: we first decouple the ODE and the PDE, then for the ODE part the equations for $\mathbf{w}$ and $\mathbf{c}$ are treated implicitly, while for the PDE part we treat the diffusion term implicitly and the nonlinear reaction term explicitly.
Given $\mathbf{v}^n$, $\mathbf{w}^n$ and $\mathbf{c}^n$ at time $t^n \in \{0, t_1, \hdots, t_{N_t}\}$,
the algebraic scheme reads as:
\begin{align*}
    &\mathbf{w}^{n+1} + \Delta t \mathbf{R}(\mathbf{v}^n,\mathbf{w}^{n+1})
    &&= \mathbf{w}^n, \\
    &\mathbf{c}^{n+1} + \Delta t \mathbf{C}(\mathbf{v}^n,\mathbf{w}^{n+1},\mathbf{c}^{n+1})
    &&= \mathbf{c}^n, \\
    &\left(\frac{\chi C_m}{\Delta t} M + A\right)\mathbf{v}^{n+1}
    &&= \frac{\chi C_m}{\Delta t} M \mathbf{v}^n
    - M \mathbf{I}_{\mathrm{ion}}(\mathbf{v}^n,\mathbf{w}^{n+1},\mathbf{c}^{n+1})
    + M \mathbf{I}_{\mathrm{app}}^{n}.
\end{align*}

Therefore, the system of ODEs governing the ionic and concentration variables is decoupled from the PDE
and solved implicitly at each timestep.
Thus, the numerical scheme requires the solution of a single linear system
associated with the parabolic equation for the transmembrane potential.
Note, in particular, that the resulting matrix $\frac{\chi C_m}{\Delta t} M + A$ is symmetric positive definite. Details on tailored solvers and preconditioners employed for each case of study will be given in Section~\ref{sec:dataset_details}.

\subsection{Surrogate Forward Architecture}
This section introduces the machine learning architecture adopted to surrogate the high-fidelity forward problem. We recall that in multi-query regimes, such as inverse problems, repeated evaluations of the forward model with slightly perturbed inputs can lead to overall prohibitive computational cost.
This motivates the use of accurate and computationally cheap surrogate models based on machine learning, such as the approach proposed in this work.
\begin{figure}[t]
    \centering
    \includegraphics[width=\textwidth,trim={2.4cm 0 2.2cm 0}, clip]{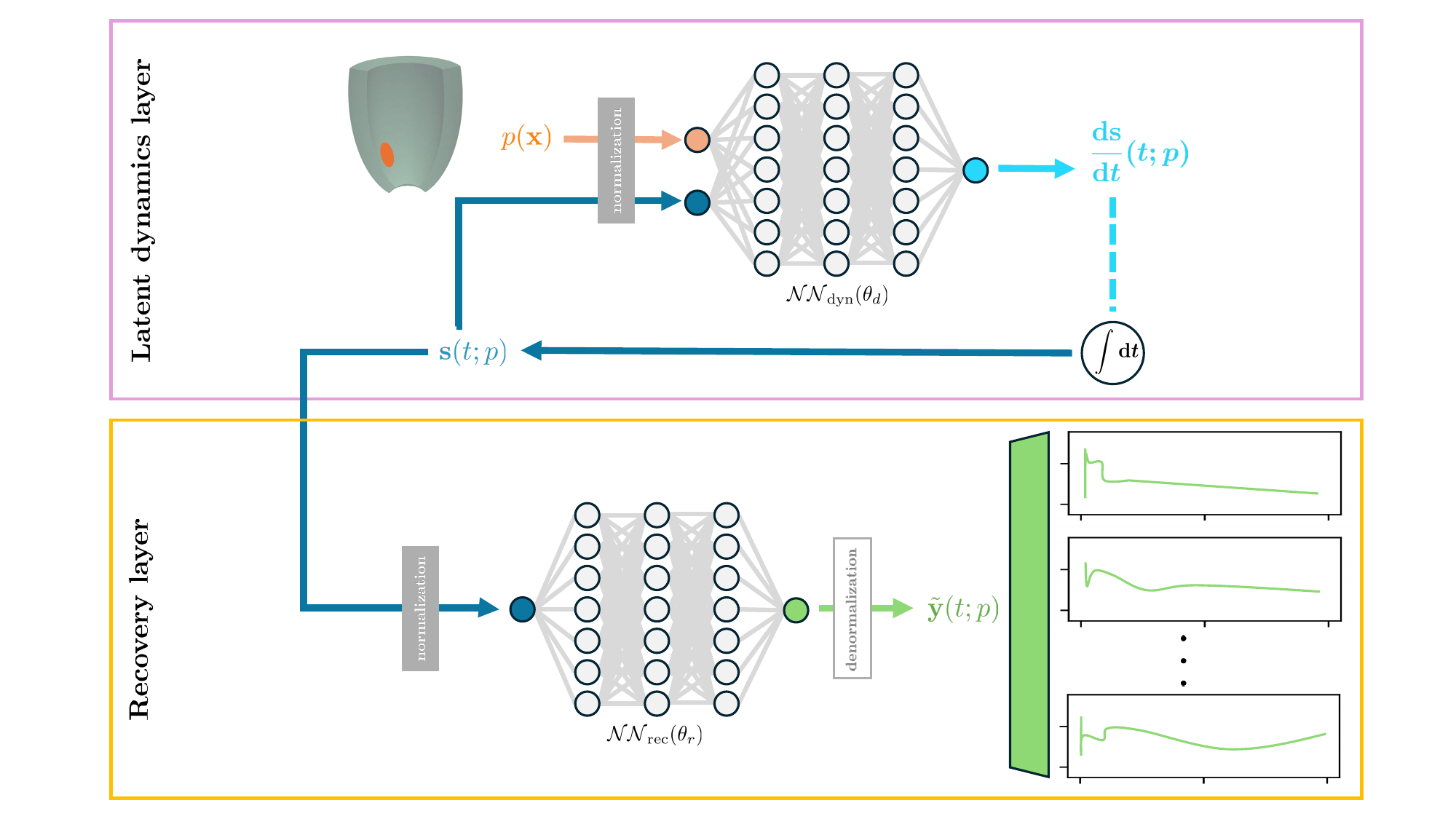}
    \caption{Forward architecture, inspired by LDNets \cite{reg24}. The latent dynamics layer is constituted by a neural ODE modeling the dependency of the latent state variables on the input parameters. In the most general case, the recovery layer can also take spatial coordinates as additional input. In the recovery layer, the latent variables are processed via a recovery feed-forward neural network to retrieve the evolution of ECG signals located at different sites. Weights and biases of the two networks are trained simultaneously by minimizing the mean squared mismatch between observation signals and reconstructions.}
    \label{fig:fwd_arch}
\end{figure}

Let $\Omega \subset \mathbb{R}^d$, $d\in\{2,3\}$, denote a bounded space domain and let $[0,T]$ be a finite time interval. Given a set of parameters $p \in \mathcal{P}$
representing spatially localized initial conditions or tissue properties (e.g.\ stimulus location or ischemic region geometry), we can define the operator,

\[
\mathcal{G}: \mathcal{P} \longrightarrow \mathcal{M}
\]
where
\[
\mathcal{G}(p) = \text{pECG}(t;p) \in \mathcal{M} \approx \mathbb{R}^{N_{\text{leads}} \times N_t},
\]
through the computation of the transmembrane potential $v(t,\mathbf{x};p)$.
Our goal is to construct a fast and reliable approximation of the forward map that preserves the essential input--output structure, rather than aiming at high-fidelity accuracy. Specifically, we search for a surrogate operator
\[
\mathcal{S} \approx \mathcal{G}, \,\mathcal{S}: \mathcal{P} \longrightarrow \mathcal{M}
\]
able to reconstruct the mapping between parameters $p$ (e.g.\ initial stimulus location or ischemic region descriptors) and the corresponding pseudo-ECG signals.

The surrogate is trained in an offline phase using a limited dataset of high-fidelity simulations,
\begin{equation}
\mathcal{D}_{\text{train}} = \{(p_i, \mathcal{G}(p_i))\}_{i=1}^{N_{\text{train}}},
\end{equation}
and subsequently deployed in an online phase to enable rapid many-query evaluations.

We adopt an architecture inspired by Latent Dynamics Networks (LDNets)~\cite{reg19, reg24} and related model and operator learning approaches~\cite{azi24, lu21}, including Recurring Neural Operators (RNOs)~\cite{liu23}, which can be interpreted as a particular case of a latent neural operator producing $0$-dimensional time series. We note that similar architectures have been previously successfully employed in biomathematics applications~\cite{zia25, zia26}, but not for inverse problems.
A schematic representation of this architecture can be found in Figure~\ref{fig:fwd_arch}.
In this setting, to construct the surrogate operator $\mathcal{S}$, we first define a set of latent state variables
\[
\mathbf{s}(t;p) \in \mathbb{R}^{n_s},
\]
whose evolution is governed by a dynamical system,
\begin{equation}
\frac{d\mathbf{s}}{dt}(t;p) = \mathcal{NN}_{\text{dyn}}(\mathbf{s}(t;p), p),
\end{equation}
where $\mathcal{NN}_{\text{dyn}}$ is a fully connected neural network to be learned.
At the discrete level, time integration is performed numerically by forward Euler method, yielding a discrete latent trajectory $\{\mathbf{s}(t_n;p)\}_{n=1}^{N_t}$.

A second neural network, called  reconstruction network,
\begin{equation}
\mathcal{NN}_{\text{rec}}: \mathbb{R}^{n_s} \longrightarrow \mathbb{R}^{N_{\text{leads}}},
\end{equation}
maps the latent state to the observable space, producing the surrogate prediction
\[
\widetilde{\text{pECG}}(t;p) = \mathcal{NN}_{\text{rec}}(\mathbf{s}(t;p)).
\]
The set of trainable variables of the two neural networks involved in the architecture $(\theta_d, \theta_r)$ are trained by minimizing the normalized MSE between ECG ground truth and reconstructions in the training set:
\begin{equation}
\begin{split}
\displaystyle \argmin_{\theta_d, \; \theta_r} \mathcal{L}(\mathcal{S}(p),\{ \text{pECG}^k\}_{k})
&= \mathrm{MSE}\big(\{\widetilde{\mathrm{pECG}^k}\}_k, 
\{\text{pECG}^k\}_{k}\big)\\
&= \frac{1}{N_{\text{leads}} N_t N_{\mathrm{samples}}}
\sum_{k=1}^{N_{\mathrm{samples}}}
\sum_{i=1}^{N_{\text{leads}}}
\sum_{j=1}^{N_t}
\left|
\text{pECG}^k_i(t_j)
-
\widetilde{\text{pECG}}^k_i(t_j;p)
\right|^2.
\end{split}
\label{eq:trainloss}
\end{equation}

In this case, spatial coordinates, corresponding to distinct leads positions, are dealt in different vectors of the output field.

\subsection{Inverse Problem: formulation and methods}\label{sec:inv_prob_form}

The inverse problem consists in reconstructing unknown parameters $p \in \mathcal{P}$ from a target set of observed pseudo-ECG measurements, i.e.
\[
\{\text{pECG}^k\}_{k \leq N_{\mathrm{samples}}} \in \mathcal{M}.
\]
This problem can be viewed as the inversion of the forward operator $\mathcal{F}$, which is generally ill-posed due to nonlinearity and a limited amount of observations.

Using the surrogate forward model, we can address the inverse problem through a standard iterative procedure, i.e. by formulating a minimization problem of a discrepancy measure for each sample to reconstruct \cite{bab04}:
\begin{equation}
\label{eq:minprob}
\min_{p \in \mathcal{P}} \;
\mathcal{J}(p)
:= \mathcal{L}\bigl(\mathcal{S}(p), \text{pECG}\bigr),
\end{equation}
where $\mathcal{L}$ is a data-misfit term for each fixed pECG datum (identified by $^{\bar{k}}$ index), which throughout this work is defined as a mean squared error loss,
\begin{equation}
\begin{split}
\mathcal{L}(\mathcal{S}(p), \text{pECG}^{\bar{k}})
&= \mathrm{MSE}\big(\text{pECG}^{\bar{k}}, \widetilde{\mathrm{pECG}^{\bar{k}}}\big)\\
&= \frac{1}{N_{\text{leads}} N_t}
\sum_{i=1}^{N_{\text{leads}}}
\sum_{j=1}^{N_t}
\left| \text{pECG}_i^{\bar{k}}(t_j) - \widetilde{\text{pECG}}_i^{\bar{k}}(t_j;p) \right|^2.
\end{split}
\label{eq:loss}
\end{equation}

The initial guess for the optimization problem~\eqref{eq:minprob} is chosen through a multiple-shooting strategy (detailed in Section~\ref{sec:guess}) and then solved using a gradient-based method, combining Adam iterations with second order quasi-Newton updates (BFGS). At each iteration, the surrogate forward model is evaluated to generate $\mathcal{S}(p)$, enabling rapid exploration of the parameter space.

For the 3-dimensional case considered, employing a non-convex hollow ellipsoid, we adopted a projected gradient descent step \cite{noc06} at each iteration, in order to avoid the algorithm to generate non-physical points outside the domain.

\section{Numerical experiments}\label{sec:num_exp}

\subsection{Dataset Details}\label{sec:dataset_details}

\begin{table}[ht]
\centering
\resizebox{\textwidth}{!}{
\begin{tblr}{
  width = \textwidth,
  colspec = {l l l},
  row{1} = {font=\bfseries},
}
\hline
Parameter & 2d configurations & 3d configuration \\
\hline

Domain geometry 
& Rectangle $(a \times b)$ 
& {Ellipsoidal parameterization \\ ($a$,$b$,$c$ in cm, $r$ constant, $\theta$,$\phi$ in rad)} \\

Domain dimensions 
& $a=5.12\,\mathrm{cm}$, $b=0.96\,\mathrm{cm}$ 
& {$a\in[2.2,3.3]$, \\ $b\in[2.2,3.3]$, \\ $c\in[5.9,6.4]$ \\ $r\in[0,1]$ \\ $\theta\in[-3\pi/8,\,\pi/8]$ \\ $\phi\in[-3\pi/2,\,\pi/2]$ } \\

Spatial discretization 
& $256 \times 48$ elements, $h=10^{-2}\,\mathrm{cm}$ 
& $(n_i,n_j,n_k)=(512,384,64)$ \\

Surface to volume ratio
& $\chi = 10^3 \mathrm{cm}^{-1}$
&  $\chi = 10^3 \mathrm{cm}^{-1}$ \\

Membrane capacitance 
& $C_m = 1\, \mu \mathrm{F}/\mathrm{cm}^2$ 
& $C_m = 1\, \mu \mathrm{F}/\mathrm{cm}^2$ \\

Ionic model 
& Ten Tusscher~\cite{ten06} 
& Ten Tusscher~\cite{ten06} \\

Ischemic radius 
& $r \sim \mathcal{U}(1.80\text,3.33)\,[\text{cm}]$ 
& -- \\

Modified ionic model (ischemia) 
& ATP-sensitive $K^+$ current \cite{dut17} 
& -- \\

Extracellular potassium (ischemia)
& $K_o = 12.0\,\mathrm{mM}$ 
& -- \\

Sodium conductance (ischemia)
& $g_{\mathrm{Na}} = 10.3866\,\mathrm{mS}/p\mathrm{F}$ 
& -- \\

L-type calcium conductance (ischemia)
& $g_{\mathrm{CaL}} = 2.78\times10^{-5}\,\mathrm{mS}/p\mathrm{F}$ 
& -- \\

Stimulus current 
& $I_{\mathrm{app}}=350\,p\mathrm{A}/\mathrm{cm}^2$ 
& $I_{\mathrm{app}}=350\,p\mathrm{A}/\mathrm{cm}^3$ \\

Time discretization 
& $\Delta t=0.05\,\mathrm{ms}$, $T_{\mathrm{end}}=400\,\mathrm{ms}$ 
& $T_{\mathrm{stim}}=1\,\mathrm{ms}$ \\

Dataset size 
& $100/100/200$ (train/val/test) 
& $300/100/100$ (train/val/test)  \\
\hline
\end{tblr}
}
\caption{Geometric, physical, and numerical parameters used for dataset generation.}
\label{tab:dataset_params}
\end{table}

Although the framework is general enough to be extended to deal with input parameters varying in functional spaces, we focus on finite-dimensional parameter spaces considering the applications of this work. In particular, using the notation proposed in Section~\ref{sec:inv_prob_form}, we will consider the following test cases:
\begin{itemize}
\item \textbf{Stimulus localization} (2d/3d): $p \in \mathbb{R}^{2,3}$ identifies the cartesian coordinates of the center of an applied stimulus.
\item \textbf{Ischemic region localization} (2d): $p \in \mathbb{R}^2$ identifies the cartesian coordinates of the center of a circular ischemic zone with fixed radius.
\item \textbf{Ischemic region with variable radius} (2d): $p \in \mathbb{R}^3$ encodes the cartesian coordinates of the center and radius of an (idealized) circular ischemic region.
\end{itemize}

The experimental dataset is partitioned into four computational scenarios, each modeling a specific cardiac pathology: 2d ectopic activation, 2d ischemia with fixed spatial parameters, 2d ischemia with variable radius, and a 3d ectopic stimulus model.

All data are generated by solving the Monodomain equations~\eqref{eq:monodomain} coupled with the Ten Tusscher ionic model, using the geometric, physical, and numerical parameters summarized in Table~\ref{tab:dataset_params}, while the pseudo-ECGs have been computed through model~\eqref{eq:pecg}.
For the 2d configurations, simulations are performed on a rectangular domain $\Omega \subset \mathbb{R}^2$ with dimensions $(a,b)$ and uniform spatial discretization. Electrical propagation is governed by anisotropic conductivity tensors, whose longitudinal and transverse intra- and extracellular coefficients are set to
\[
\sigma_{l}^e=2\times10^{-3}, \quad
\sigma_{t}^e=1.3514\times10^{-3}, \quad
\sigma_{l}^i=3\times10^{-3}, \quad
\sigma_{t}^i=3.1525\times10^{-4}
\;\; p\mathrm{F}/\mathrm{cm}^2.
\]

An external stimulus current $I_{\mathrm{app}}$ is applied to initiate activation, and the system is integrated in time using a fixed time step up to a prescribed final time.
In the ischemic and ischemic-radius configurations, tissue heterogeneity is introduced through a modified Ten Tusscher ionic model including an ATP-sensitive potassium current \cite{dut17}. In this setting, the extracellular potassium concentration $K_o$, sodium conductance $g_{\mathrm{Na}}$, and L-type calcium conductance $g_{\mathrm{CaL}}$ are modified as reported in Table~\ref{tab:dataset_params}. In the ischemic-radius case, the spatial extent of the ischemic region is modeled by a circular inclusion whose radius is treated as a random parameter sampled uniformly in the interval $[1.80,\,3.33]\, (\text{cm})$.
Each 2d configuration is split into training, validation, and test sets of fixed size, see Table \ref{tab:dataset_params}.
Each high-fidelity solution has been computed with a MATLAB code, where the  linear system resulting from the FEM discretization~\eqref{eq:mono-matrix} is solved with standard  \textit{backslash}. Note that a single solution run on an M1-chip equipped laptop required about 8 minutes.
The 3d dataset is generated on an ellipsoidal geometry $\Omega \subset \mathbb{R}^3$ defined through the parametric mapping:
\begin{equation}
\label{eq:ellipsoid_param}
\begin{cases}
   x = a(r)\cos\theta\cos\phi, \\
   y = b(r)\cos\theta\sin\phi, \\
   z = c(r)\sin\phi,
\end{cases}
\end{equation}
where $r\in[0,1]$, $\theta\in[-3\pi/8,\,\pi/8]$, and $\phi\in[-3\pi/2,\,\pi/2]$. 
The functions $a(r)$, $b(r)$, and $c(r)$ interpolate between the semi-axis bounds $a\in[a_1,a_2]$, $b\in[b_1,b_2]$, and $c\in[c_1,c_2]$ reported in Table~\ref{tab:dataset_params}. 
The parameterized domain is discretized on a structured grid in the $(r,\theta,\phi)$ coordinates. Electrical activation is initiated by a short-duration stimulus current, and electrophysiological dynamics are again modeled using the Ten Tusscher ionic model. Training, validation, and test samples are obtained by uniformly sampling observation points on the 3d geometry.
Each high-fidelity solution has been solved using a PETSc~\cite{bal23} implementation. We notice that in this case the high number of degrees of freedom of the problem required employing an iterative solver. In particular, since the iteration matrix derived from the IMEX scheme described in Section~\ref{sec:num_schemes} is symmetric definite positive, we employed the conjugate method (CG) preconditioned with the algebraic multigrid (AMG) provided by HYPRE~\cite{fal06}, wrapped in PETSc.
From a computational standpoint, each high-fidelity simulation required around 20 minutes on a single node of CINECA HPC supercomputer LEONARDO~\cite{tur24}, equipped with 4 NVIDIA A100 GPUs each with 32GB memory.

\begin{figure}[t]
    \centering
    \includegraphics[trim={5cm 0 5cm 0},clip,width=0.7\linewidth]{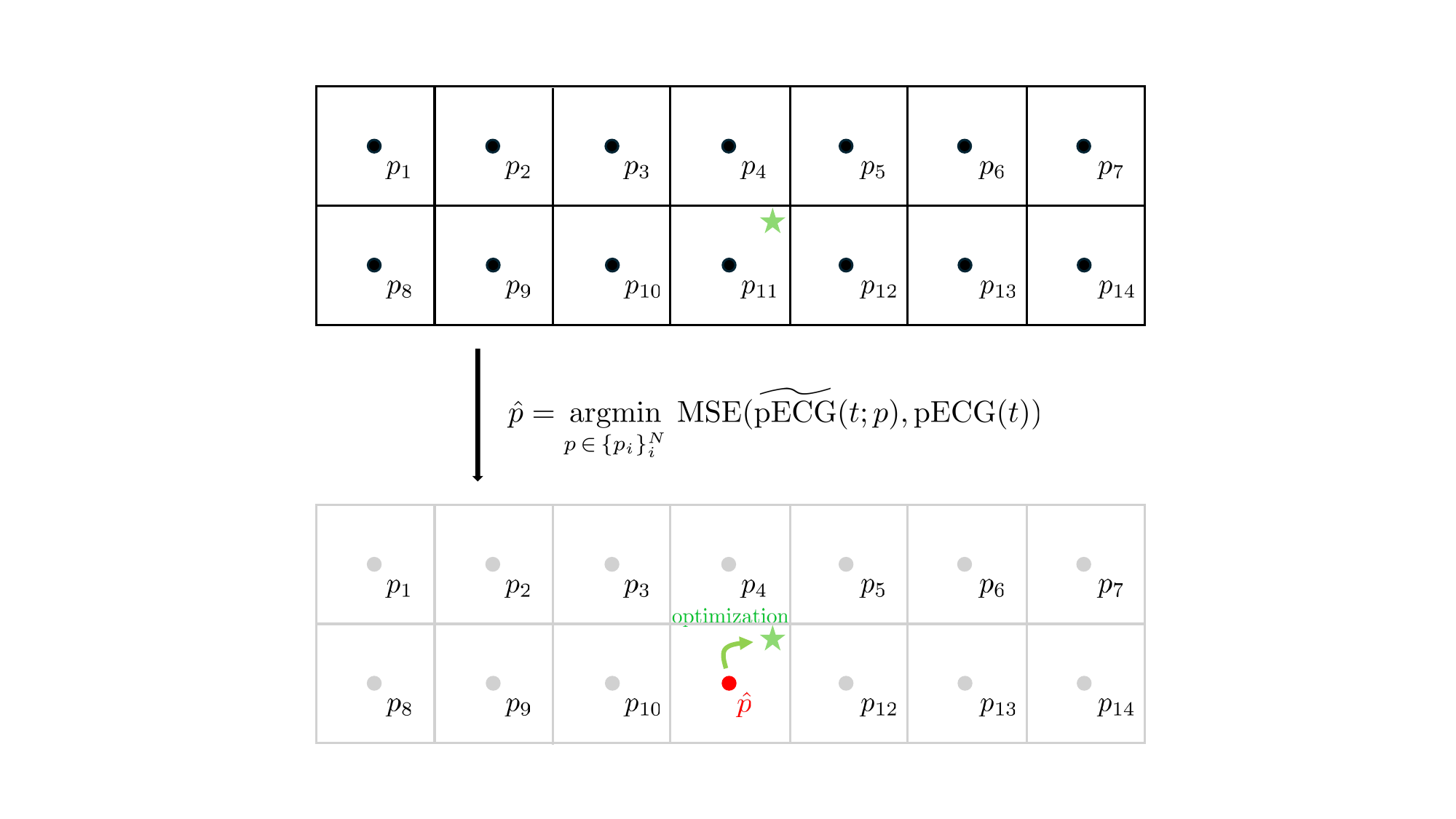}
    \caption{Multi-start strategy for the 2d domain. The starred point is the target: the initial guess of the inverse problem loop is chosen as the center of the subdomains which generates the pseudo-ECG closest to the observed one.}
    \label{fig:dd2d}
\end{figure}

\begin{figure}[t]
    \centering
    \includegraphics[width=0.8\textwidth]{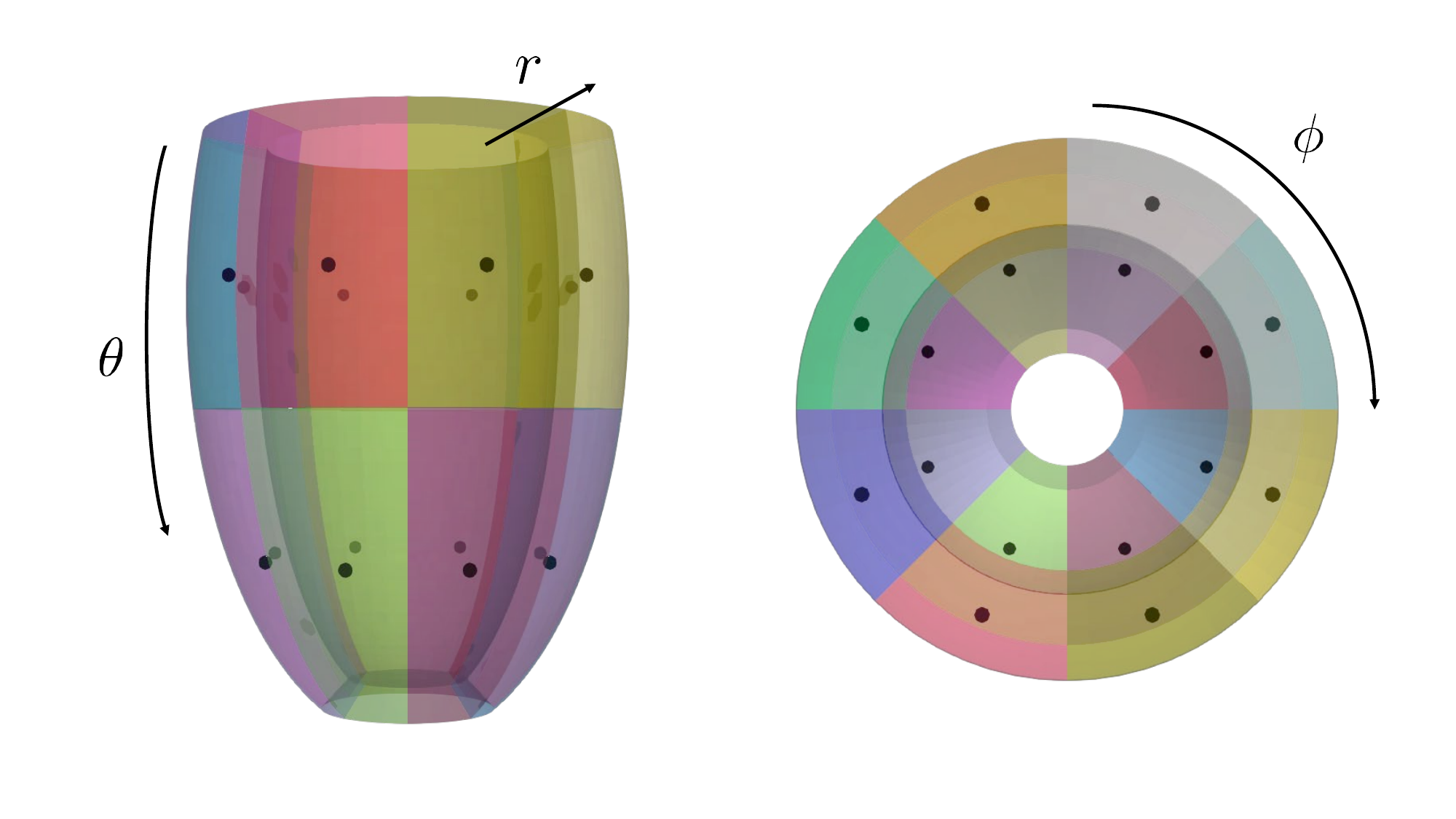}
    \caption{Representation of repartition in subdomains of the 3d domain for the multi-start strategy. The partitioning has been performed using spherical coordinates.}
    \label{fig:dd3d}
\end{figure}

\subsection{Initial guess for the inverse problem solution} \label{sec:guess}
 
Training is carried out using a three-stage Adam optimization algorithm with progressively decreasing learning rates, followed by a second order BFGS refinement phase.

To initialize the inverse optimization procedure, we adopt a multi-start strategy based on a discrete set of candidate points distributed over the computational domain (see Figure~\ref{fig:dd2d} for the 2d case and Figure~\ref{fig:dd3d} for the 3d case). 

Let $\{p_i\}_{i=1}^N$ denote a predefined set of candidate locations, which may be chosen by uniformly subdividing the domain and picking the center of each of the subdomains. 
Regarding the 2d case, for each $p_i$, we evaluate the discrepancy between the corresponding simulated pseudo-ECG signal and the observed one, i.e.
\[
\displaystyle \hat{p} 
= \argmin_{p \, \in \, \{p_i\}_i^N} 
\mathrm{MSE}\big(\widetilde{\mathrm{pECG}}(t,p), 
\mathrm{pECG}(t)\big),
\]

and we select as initial guess the point minimizing the data misfit.
This discrete screening step provides a physically meaningful prior located in the region of highest similarity with the measured signal. 
Therefore, starting from $\hat{p}$ we perform a gradient-based optimization in the parameter space to refine the estimate and recover the target location.

For the 3d case, we observed a significant impact of identifiability issues, i.e., different initial stimuli may generate very similar pseudo-ECG signals, leading the above strategy to fail in more than $20\%$ of the test cases. 
To address this limitation while maintaining a reasonable computational cost for potential clinical applications, we adopt a more robust approach: (i) we perform some Adam iterations of the inverse optimization starting from each candidate point $\{p_i\}_i$ in the domain partition, and (ii) we conclude the inverse optimization starting from the point which minimizes the discrepancy with the observed pseudo-ECG after step (i). 
As discussed in Section~\ref{sec:2d}, this strategy is an acceptable trade-off between computational efficiency and reconstruction accuracy in the presence in non-convex-optimization scenarios.

\subsection{ECG-initial stimulus (2d/3d)}\label{sec:2d}

\begin{figure}[t]
    \centering
    \begin{subfigure}{0.46\textwidth}
        \centering
        \includegraphics[width=\linewidth,trim={0cm 0cm 1cm 0cm},clip]
        {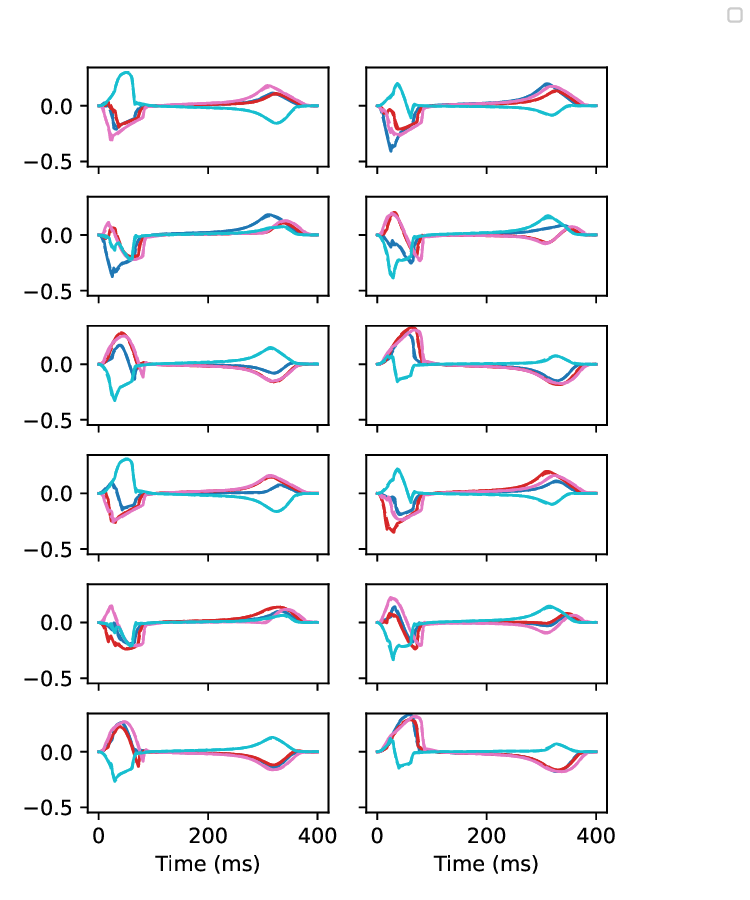}
        \caption{Forward problem. Each plot represents potential vs time for a lead for the pseudo-ECG solution.}
    \end{subfigure}
    \hfill
    \begin{subfigure}{0.482\textwidth}
        \centering
        
        \begin{subfigure}{\textwidth}
            \centering
            \includegraphics[width=\linewidth,
            trim={1.2cm 4.2cm 1.5cm 4.8cm},clip]
            {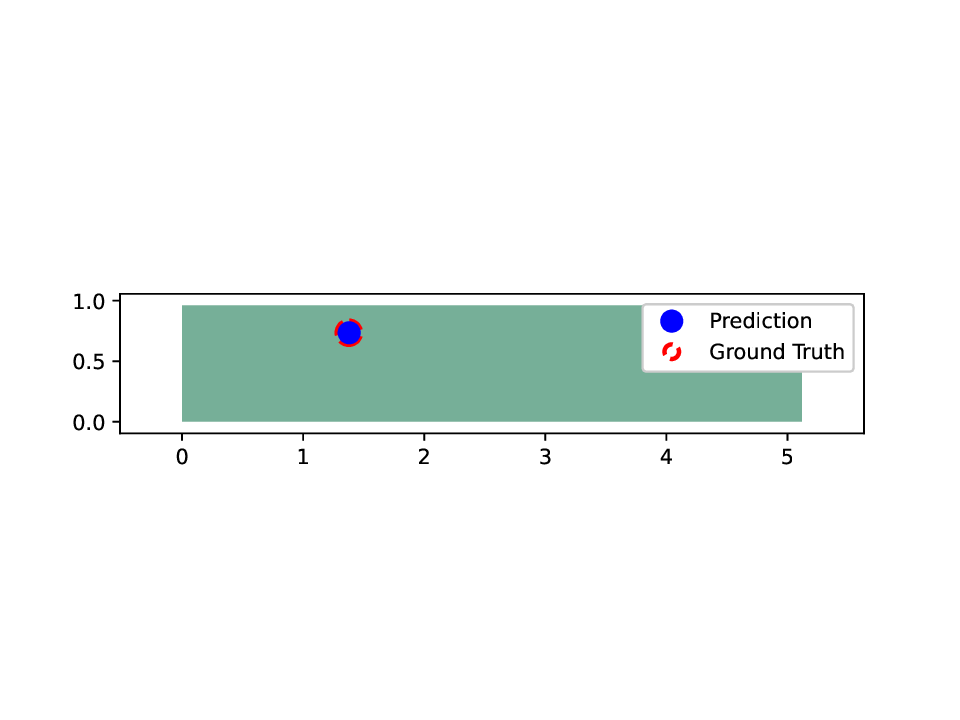}
        \end{subfigure}
        
        \vspace{0.3cm}
        
        \begin{subfigure}{\textwidth}
            \centering
            \includegraphics[width=\linewidth,
            trim={1.2cm 4.2cm 1.5cm 4.8cm},clip]
            {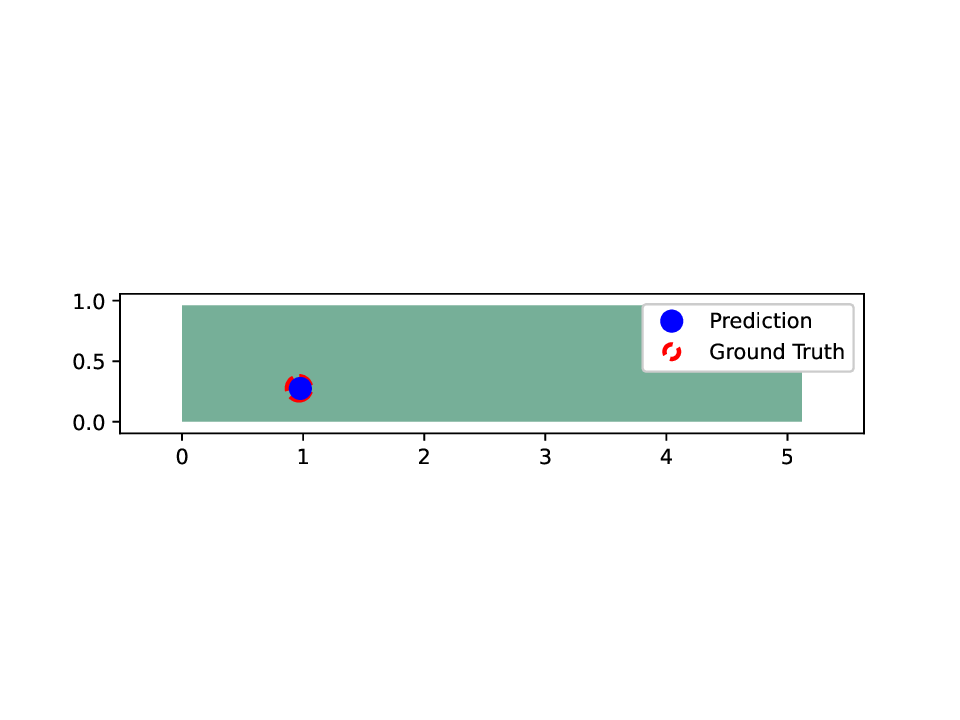}
        \end{subfigure}
        
        \vspace{0.3cm}
        
        \begin{subfigure}{\textwidth}
            \centering
            \includegraphics[width=\linewidth,
            trim={1.2cm 4.2cm 1.5cm 4.8cm},clip]
            {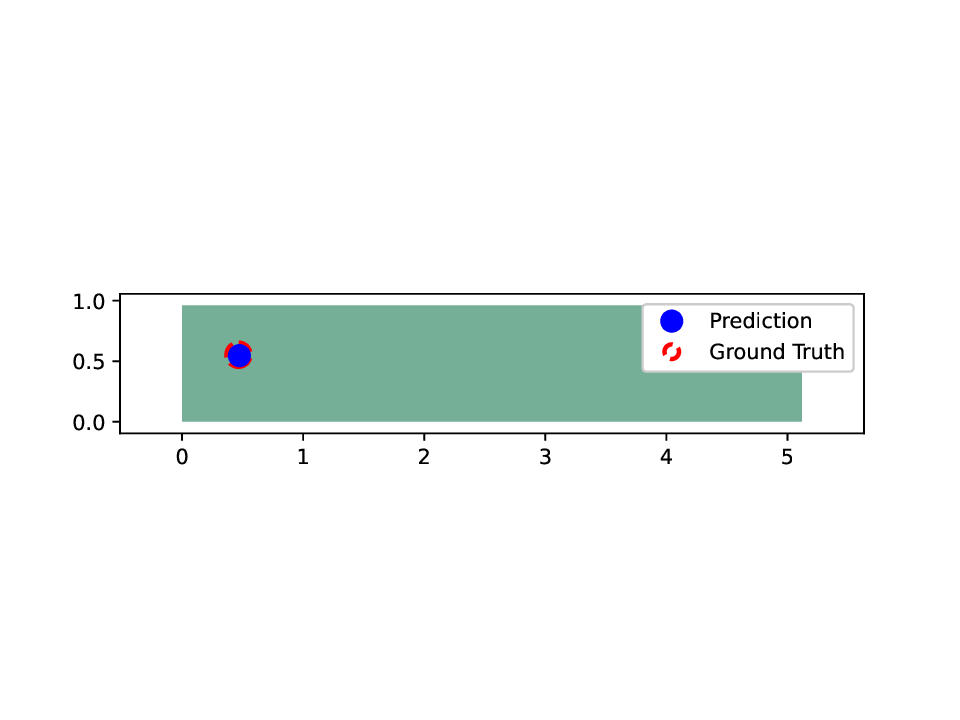}
        \end{subfigure}
        
        \vspace{0.3cm}
        
        \begin{subfigure}{\textwidth}
            \centering
            \includegraphics[width=\linewidth,
            trim={1.2cm 4.2cm 1.5cm 4.8cm},clip]
            {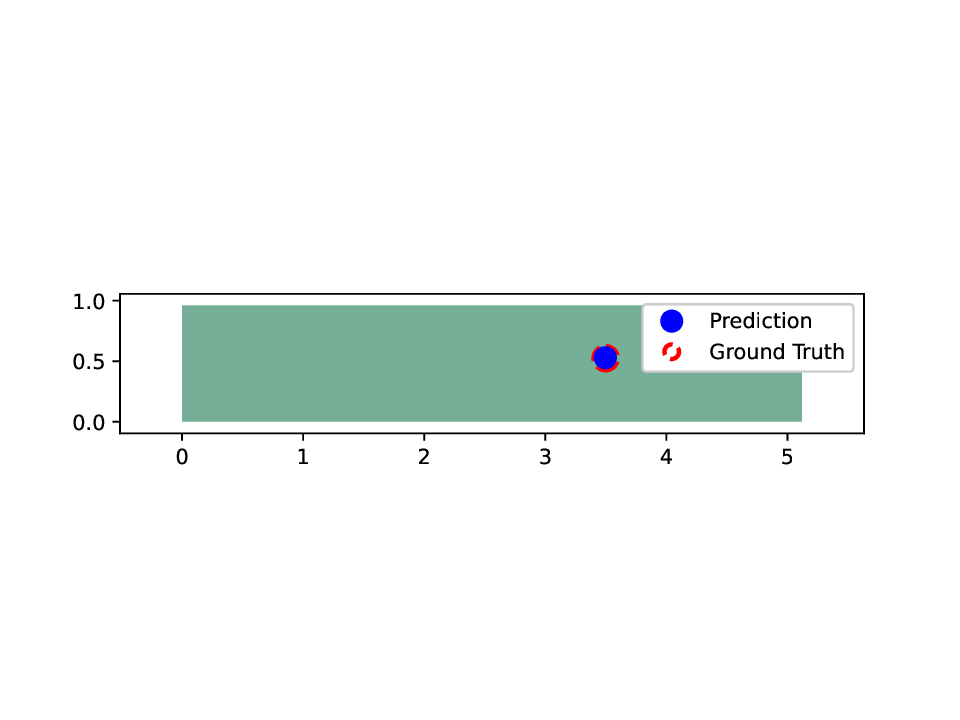}
        \end{subfigure}

        \vspace{0.65cm}
        
        \caption{Inverse problem reconstructions corresponding to the forward solutions shown on the left.}
        \label{fig:vanilla2d_inverse}
    \end{subfigure}
    
    \caption{
    2d initial stimulus test case. 
    \textbf{Left (color online):} Forward problem solutions (pseudo-ECGs generated by the ROM for different parameter instances). Different colors indicate different cases. Dashed lines indicates ground truth solution, while continuous lines indicate ML reconstructions. Results show good correspondence between ground truth and neural network approximation.
    \textbf{Right (color online):} Corresponding inverse reconstructions obtained from those forward-generated signals. 
    Each inverse solution is computed using the pseudo-ECG displayed on the left as observed data.
    }
    \label{fig:vanilla2d_combined}
\end{figure}

\begin{figure}
    \centering
    
    \begin{subfigure}{\linewidth}
        \centering
        \includegraphics[width=\linewidth]{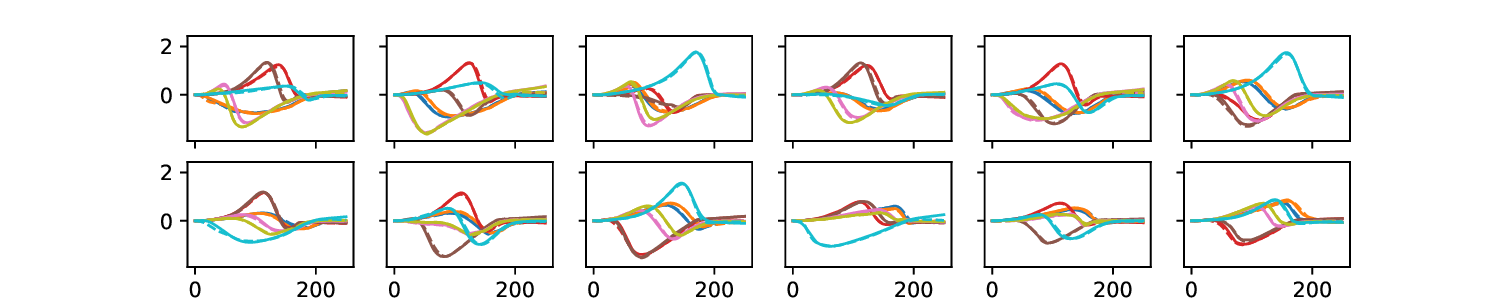}
        \caption{}
    \end{subfigure}
    
    \vspace{0.5cm}
    
    \begin{subfigure}{\linewidth}
        \centering
        \includegraphics[width=\linewidth]{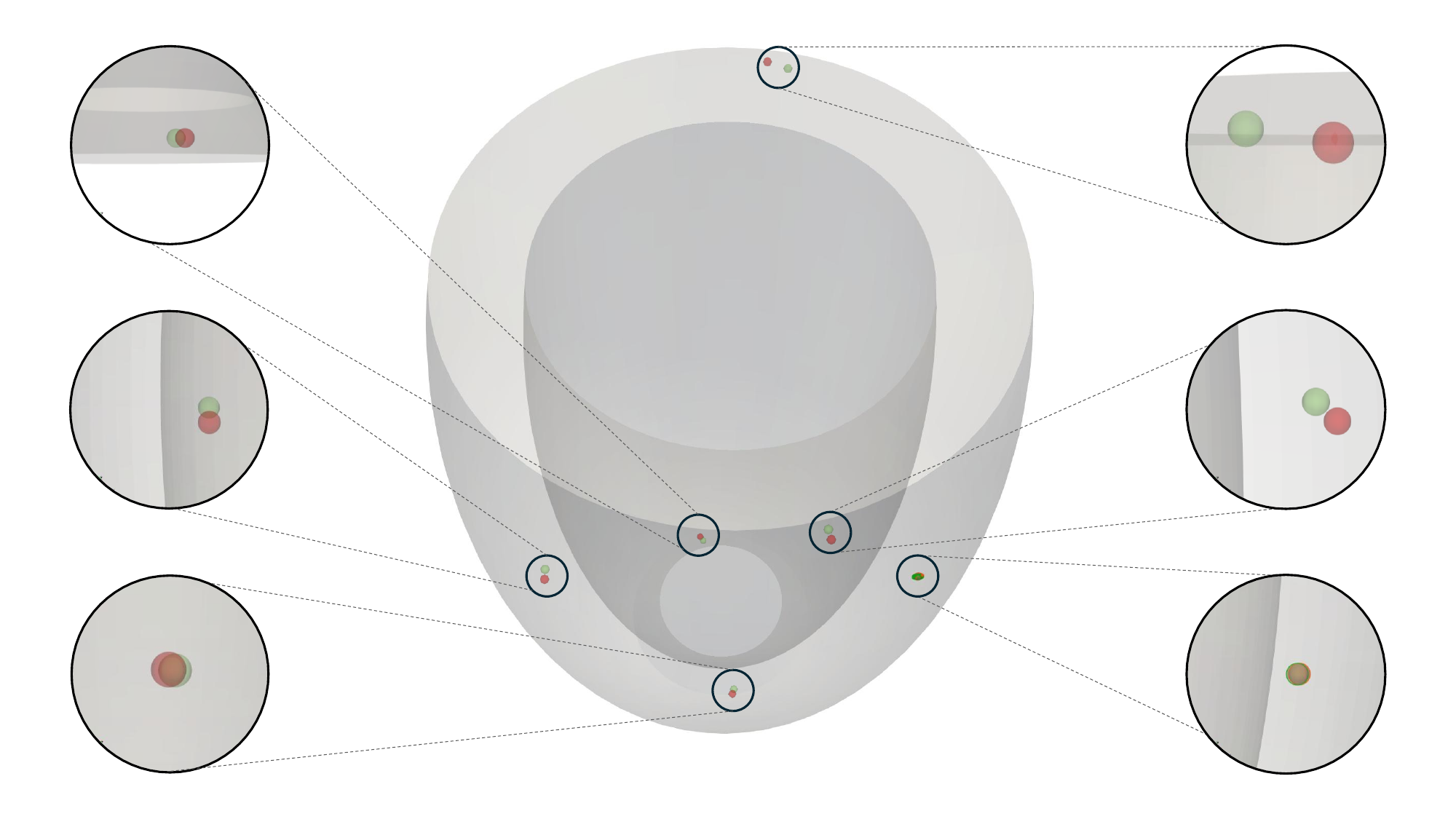}
        \caption{}
    \end{subfigure}
    
    \caption{3d initial stimulus test case. 
    \textbf{Top (color online):} Forward problem solutions (pseudo-ECGs generated by the ROM for different parameter instances). Different colors indicate different cases. Dashed lines indicates ground truth solution, while continuous lines indicate ML reconstructions. Results show good correspondence between ground truth and neural network approximation.
    \textbf{Bottom (color online):} Corresponding inverse reconstructions obtained from those forward-generated signals (green is the ground truth while red is the corresponding reconstruction). 
    Each inverse solution is computed using the pseudo-ECG displayed on the left as observed data.}
    \label{fig:vanilla3dcombined}
\end{figure}

\begin{figure}[t]
    \centering
    
    \begin{subfigure}{0.4\linewidth}
        \centering
        \includegraphics[width=\linewidth]{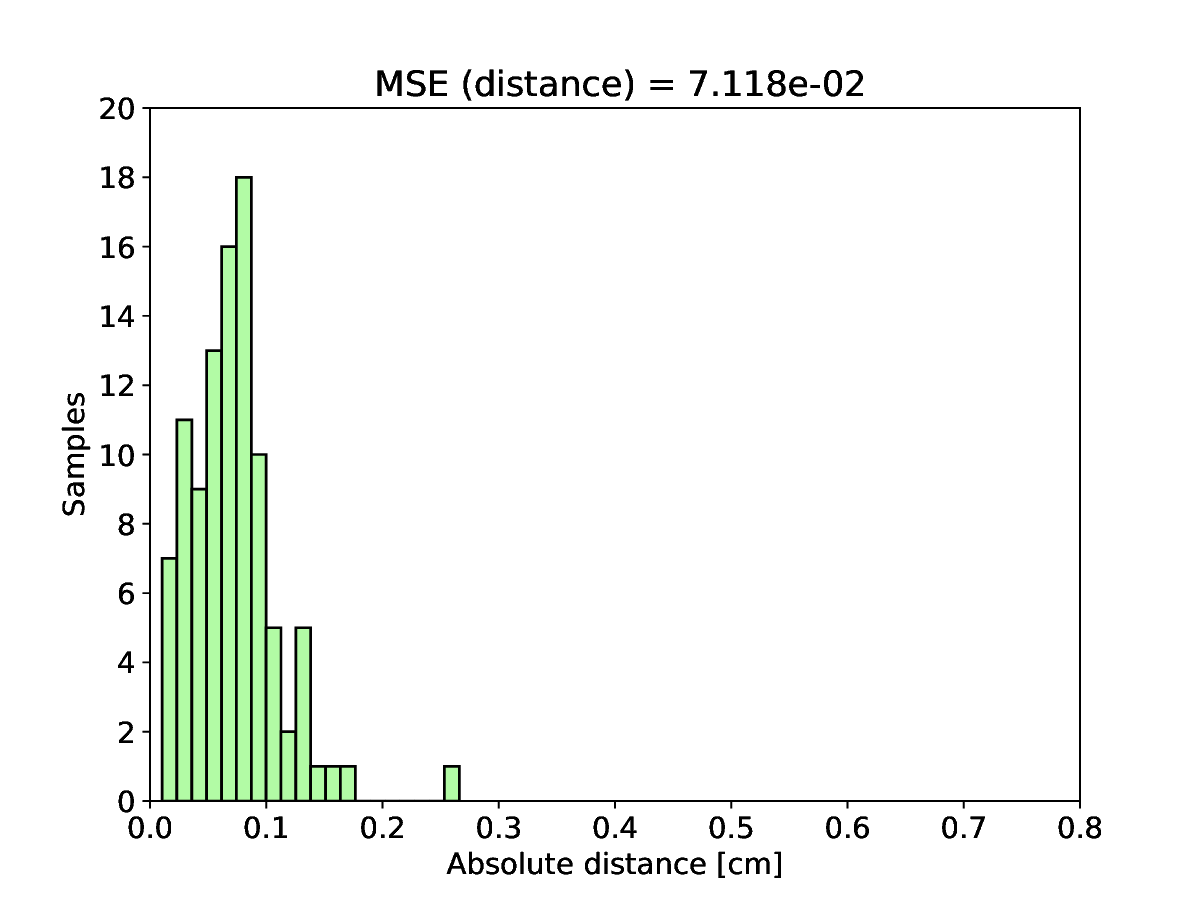}
        \caption{}
    \end{subfigure}
    \qquad
    \begin{subfigure}{0.4\linewidth}
        \centering
        \includegraphics[width=\linewidth]{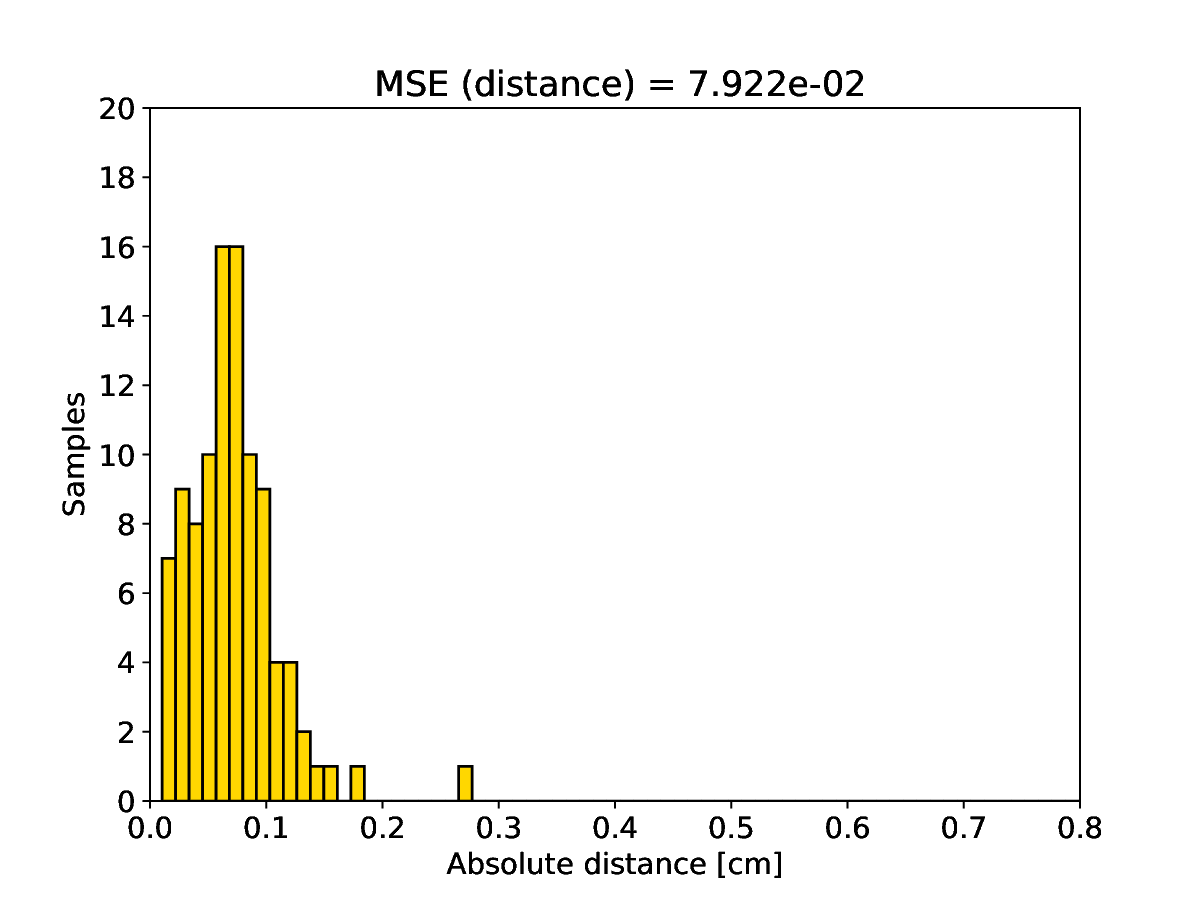}
        \caption{}
    \end{subfigure}
    \vspace{0.1cm}
    
    \begin{subfigure}{0.4\linewidth}
        \centering
        \includegraphics[width=\linewidth]{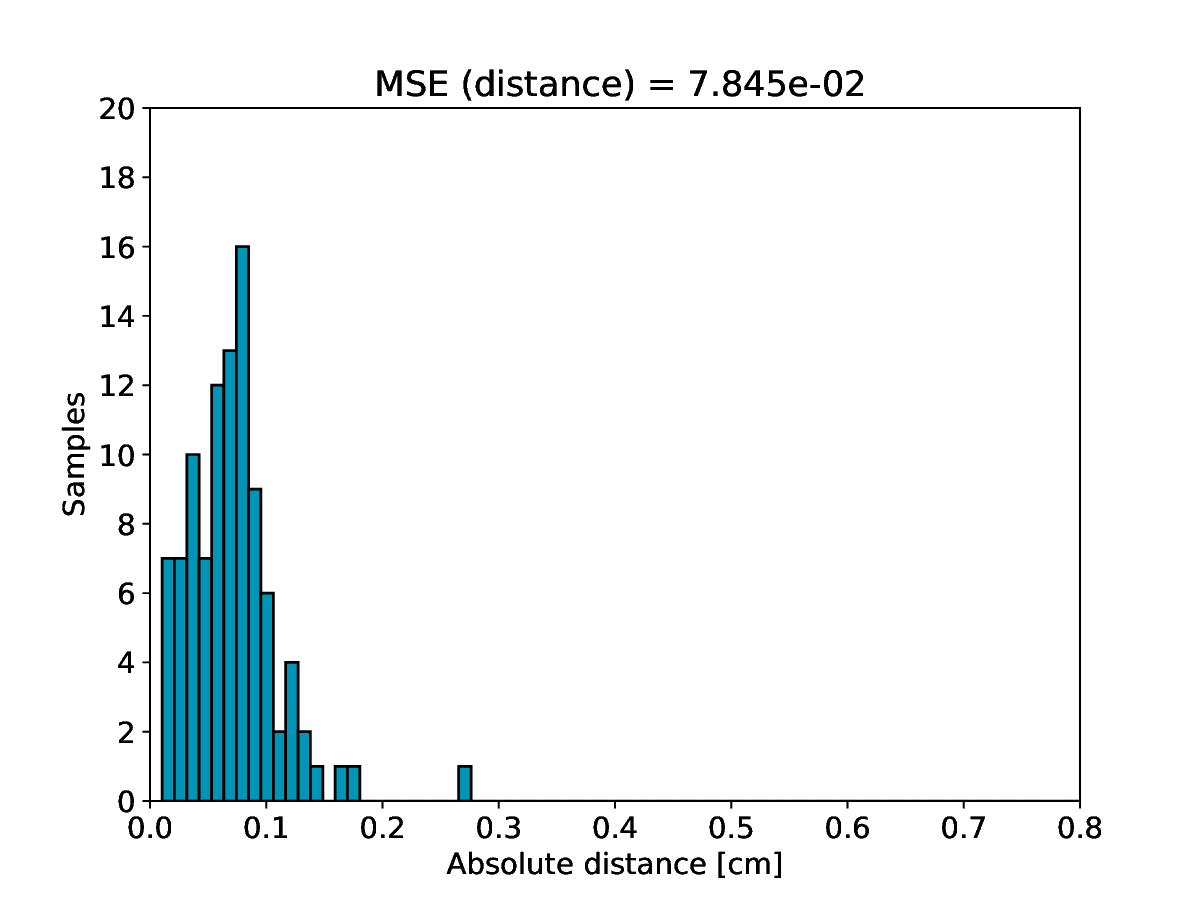}
        \caption{}
    \end{subfigure}
    \qquad
    \begin{subfigure}{0.4\linewidth}
        \centering
        \includegraphics[width=\linewidth]{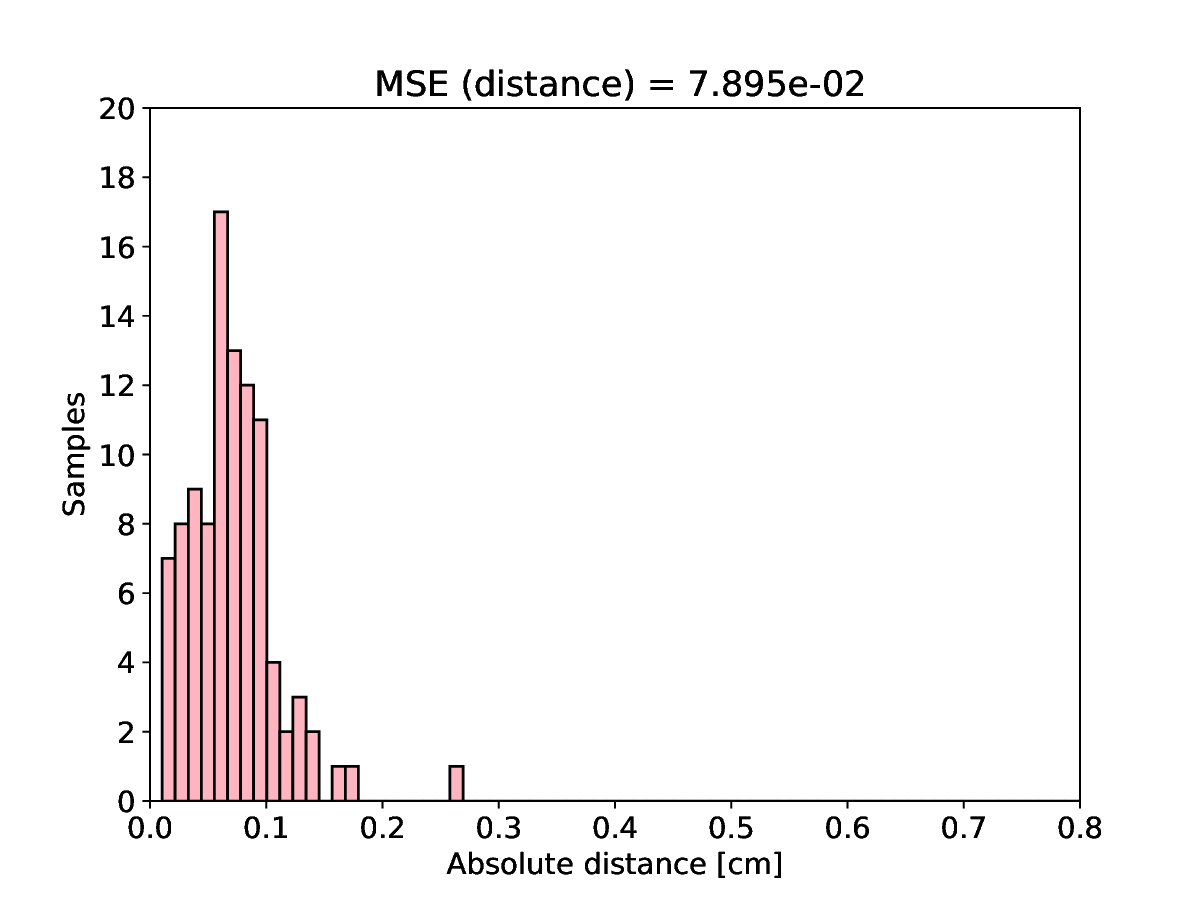}
        \caption{}
    \end{subfigure}
    
    \caption{Mean square error of the distance between the ground truth stimulus center and the reconstruction (3d case). We consider the following number of subdivisions of the domain in the radial direction (a) 1, (b) 2, (c) 4 and (d) 6. }
    \label{fig:msedists}
\end{figure}

In Table~\ref{tab:arch_params_2d} we can find the experimental setup for the 2d case.
Results reported in Tables~\ref{tab:latent_states_2d}–\ref{tab:training_stats_2d} assess both the accuracy of the forward latent dynamics model and the effectiveness of the inverse optimization procedure used to localize the initial activation site.
In particular, Table~\ref{tab:latent_states_2d} highlights how the forward model provides accurate pseudo-ECG reconstructions across a wide range of latent dimensions, with normalized RMSE values consistently of order $10^{-3}$ and Pearson dissimilarity remaining below $4\times10^{-3}$.
The best compromise between accuracy and computational cost is obtained for $12$–$16$ latent states. 
Increasing the latent dimension beyond this range does not lead to systematic improvements and may result in a degradation of generalization performance, as observed for $18$ and $20$ latent states.
However, for 24 latent states the mean value for validation loss is minimized to $1.778\times10^{-4}\pm7.777\times10^{-6}$ (notice that in Table~\ref{tab:latent_states_2d} only the smallest value obtained is reported), while memory usage and training time remain within the observed range. Hence, we choose to deal with a model with 20 latent states for surrogating the forward problem.

The inverse optimization loop for 2d stimulus localization exhibits stable and reliable convergence behavior, as summarized in Table~\ref{tab:training_stats_2d}.
The minimum achieved loss reaches values below $10^{-7}$, while the mean loss over the test dataset remains on the order of $10^{-4}$ with limited variability across samples: this aspect entails robust recovery of the initial stimulus location from pseudo-ECG measurements.
Furthermore, the inverse problem is solved with modest computational requirements, with an average optimization time of approximately $80$ seconds per instance and limited CPU memory usage, confirming the practical feasibility of the proposed approach for 2d ECG-based stimulus reconstruction.
Plots of the forward and the inverse problem solutions for the 2d case are reported in Figure~\ref{fig:vanilla2d_combined}.

For the 3d ECG–initial stimulus reconstruction task we provide details regarding the architecture in Table~\ref{tab:arch_params3d}.
Results reported in Tables~\ref{tab:latent_states_3d}-\ref{tab:training_stats3d} underline both the increased complexity of the forward approximation and the effectiveness of the inverse optimization procedure in this 3-dimensional case. 
Although generalization errors are higher in this case, Table~\ref{tab:latent_states_3d} shows that the forward model is able to accurately reproduce pseudo-ECG signals, reflecting the higher complexity of the 3d geometry and dynamics. 
Normalized RMSE values are of order $10^{-2}$ and Pearson dissimilarity remains below $1.6\times10^{-2}$. 
The best trade-off between accuracy and computational cost is achieved for $16$-$20$ latent states, where both error metrics are minimized. 
Further increasing the latent dimension does not yield significant improvements while substantially increasing memory usage and training time, which already reach considerable values due to the large-scale nature of the problem. Therefore, even in this case, the number of latent states for surrogating the forward problem is fixed at 20.

Retrieving stimulus localization in 3d remains effective despite the increased dimensionality and the presence of identifiability issues, which have been mitigated through the procedure described in Section~\ref{sec:guess}.
In Figure \ref{fig:msedists} we study the distribution of the mean squared error (MSE) across the test set as the number of radial subdivisions is varied from 1 to 6.
Our results indicate that the number of initial guesses in the transmural direction does not significantly influence the accuracy of the inverse problem solution, which consistently maintains an average error of approximately 0.0007 cm.
Conversely, increasing the number of radial subdivisions introduces outliers into the error distribution and increases the overall computational complexity of the minimization process.

As reported in Table~\ref{tab:training_stats3d}, the minimum loss achieved is of order $10^{-6}$, while the mean loss over the test dataset remains below $10^{-3}$, albeit with higher variability compared to the 2d case.

As expected, the overall computational burden is higher, but still compatible with offline or near real time clinical applications with an average optimization time of approximately $640$ seconds per instance.
Plots of the forward and the inverse problem solutions for the 3d case are summarized in Figure~\ref{fig:vanilla3dcombined}.
\begin{table}[t]
\centering
\begin{tblr}{
  width = \textwidth,
  colspec = {l c c},
  row{1} = {font=\bfseries},
}
\hline
Parameter & Forward model & Inverse model \\
\hline

Temporal step $dt$ & $6\times10^{-3}$ & -- \\

Dynamic network layers & $[4,\,8]$ & -- \\

Reconstruction network layers & $[17,\,17,\,17,\,17]$ & -- \\

Regularization parameter $\alpha$ & $4.7\times10^{-3}$ & -- \\

Number of subdomains & -- & $40 \times 20$ \\

Adam epochs (stage 1 / stage 2) & $300 / 300$ & $20 / 20$ \\

Adam learning rate (stage 1) & $1\times10^{-2}$ & $1\times10^{-2}$ \\

Adam learning rate (stage 2) & $1\times10^{-3}$ & $1\times10^{-3}$ \\

BFGS epochs & $10{,}000$ & $20$ \\
\hline
\end{tblr}
\caption{
Architecture and optimization parameters for the forward and inverse models for the 2d initial stimulus case.
Here, $dt$ denotes the temporal step size of the forward Euler scheme used in the dynamics network,
the dynamic and reconstruction network layers specify the number of neurons per hidden layer,
$\alpha$ is the regularization parameter weighting the training loss,
Adam epochs and learning rates correspond to the two-stage Adam optimization procedure,
BFGS epochs indicate the number of iterations of the quasi-Newton optimizer used for fine-tuning,
while the number of subdomains defines the spatial partitioning adopted for the optimization of the inverse problem.
}
\label{tab:arch_params_2d}
\end{table}

\begin{table}[t]
\centering
\resizebox{\textwidth}{!}{  
\begin{tblr}{
  colspec = {c c c c c c c},
  row{1} = {font=\bfseries},
}
\hline
Latent states 
& Training loss 
& Validation loss 
& Normalized RMSE 
& Pearson dissimilarity 
& CPU memory [GB] 
& Training time [h] \\
\hline
8  & $6.612\times10^{-4}$ & $2.033\times10^{-4}$ & $7.034\times10^{-3}$ & $2.857\times10^{-3}$ & 4.47 & 10.45 \\
12 & $5.855\times10^{-4}$ & $1.821\times10^{-4}$ & $6.695\times10^{-3}$ & $2.587\times10^{-3}$ & 5.38 & 9.91  \\
16 & $5.817\times10^{-4}$ & $1.888\times10^{-4}$ & $6.797\times10^{-3}$ & $2.668\times10^{-3}$ & 5.75 & 9.75  \\
18 & $7.619\times10^{-4}$ & $2.662\times10^{-4}$ & $8.018\times10^{-3}$ & $3.714\times10^{-3}$ & 5.69 & 9.76  \\
20 & $5.912\times10^{-4}$ & $2.617\times10^{-4}$ & $7.553\times10^{-3}$ & $3.295\times10^{-3}$ & 6.10 & 10.10 \\
24 & $5.400\times10^{-4}$ & $1.687\times10^{-4}$ & $6.436\times10^{-3}$ & $2.392\times10^{-3}$ & 6.61 & 10.82 \\
\hline
\end{tblr}
}
\caption{Performance metrics for the 2d initial stimulus model with varying numbers of latent states (forward problem approximation).}
\label{tab:latent_states_2d}
\end{table}

\begin{table}[t]
\centering
\begin{tblr}{
  colspec = {l l},
  row{1} = {font=\bfseries},
}
\hline
Metric & Value \\
\hline
Minimum loss & $9.25\times10^{-8}$ \\
Mean loss (test dataset) & $1.06\times10^{-4} \pm 1.06\times10^{-5}$ \\
Maximum CPU memory usage & $5.30~\mathrm{GB}$ \\
Mean optimization time per data & $82~\mathrm{s}$ \\
\hline
\end{tblr}
\caption{Summary statistics of the inverse optimization loop for the 2d initial stimulus case.}
\label{tab:training_stats_2d}
\end{table}
\begin{table}[t]
\centering
\begin{tblr}{
  width = \textwidth,
  colspec = {l c c},
  row{1} = {font=\bfseries},
}
\hline
Parameter & Forward model & Inverse model \\
\hline

Temporal step $dt$ & $10\times10^{-3}$ & -- \\

Dynamic network layers & $[24,\,20,\, 20, \,12]$ & -- \\

Reconstruction network layers & $[20,\,40,\,40,\,23]$ & -- \\

Regularization parameter $\alpha$ & $3.7\times10^{-3}$ & -- \\

Number of subdomains & -- & $4 \times 1 \times 4\; [\theta, \,r,\, \phi]$ \\

Adam epochs (stage 1/stage 2) & $10000 / 10000$ & $30 / 30$ \\

Adam learning rate (stage 1) & $1\times10^{-2}$ & $1\times10^{-2}$ \\

Adam learning rate (stage 2) & $1\times10^{-3}$ & $1\times10^{-3}$ \\

BFGS epochs & $2500$ & $20$ \\
\hline
\end{tblr}
\caption{
Architecture and optimization parameters for the forward and inverse models for the 3d initial stimulus case.
Here, $dt$ denotes the temporal step size of the forward Euler scheme used in the dynamics network,
the dynamic and reconstruction network layers specify the number of neurons per hidden layer,
$\alpha$ is the regularization parameter weighting the training loss,
Adam epochs and learning rates correspond to the two-stage Adam optimization procedure,
BFGS epochs indicate the number of iterations of the quasi-Newton optimizer used for fine-tuning,
while the number of subdomains defines the spatial partitioning adopted for the optimization of the inverse problem.
}
\label{tab:arch_params3d}
\end{table}

\begin{table}[t]
\centering
\resizebox{\textwidth}{!}{  
\begin{tblr}{
  colspec = {c c c c c c c},
  row{1} = {font=\bfseries},
}
\hline
Latent states 
& Training loss 
& Validation loss 
& Normalized RMSE 
& Pearson dissimilarity 
& CPU memory [GB] 
& Training time [h] \\
\hline
8  & $1.100\times10^{-3}$ & $8.291\times10^{-4}$ & $1.441\times10^{-2}$ & $1.051\times10^{-2}$ & 43.13 & 67.72 \\
12 & $1.759\times10^{-3}$ & $1.151\times10^{-3}$ & $1.738\times10^{-2}$ & $1.534\times10^{-2}$ & 46.55 & 65.53  \\
16 & $6.969\times10^{-4}$ & $4.023\times10^{-4}$ & $1.027\times10^{-2}$ & $5.327\times10^{-3}$ & 69.64 & 109.86  \\
20 & $6.254\times10^{-4}$ & $3.481\times10^{-4}$ & $9.471\times10^{-3}$ & $4.531\times10^{-3}$ & 67.77 & 87.56 \\
26 & $6.251\times10^{-4}$ & $3.716\times10^{-4}$ & $9.764\times10^{-3}$ & $4.816\times10^{-3}$ & 80.76 & 91.96 \\
\hline
\end{tblr}
}
\caption{Performance metrics for the 3d initial stimulus model with varying numbers of latent states (forward problem approximation).}
\label{tab:latent_states_3d}
\end{table}

\begin{table}[ht]
\centering
\begin{tblr}{
  colspec = {l l},
  row{1} = {font=\bfseries},
}
\hline
Metric & Value \\
\hline
Minimum loss & $3.24\times10^{-6}$ \\
Mean loss (test dataset) & $4.60\times10^{-4} \pm 2.50\times10^{-3}$ \\
Maximum CPU memory usage & $5.45~\mathrm{GB}$ \\
Mean optimization time per data & $640~\mathrm{s}$ \\
\hline
\end{tblr}
\caption{Summary statistics of the inverse optimization loop for the 3d initiali stimulus case.}
\label{tab:training_stats3d}
\end{table}
\

\

%

\subsection{ECG-ischemic region (2d)}

\begin{figure}[t]
    \centering
    \begin{subfigure}{0.46\textwidth}
        \centering
        \includegraphics[width=\linewidth,trim={0cm 0cm 1cm 0cm},clip]
        {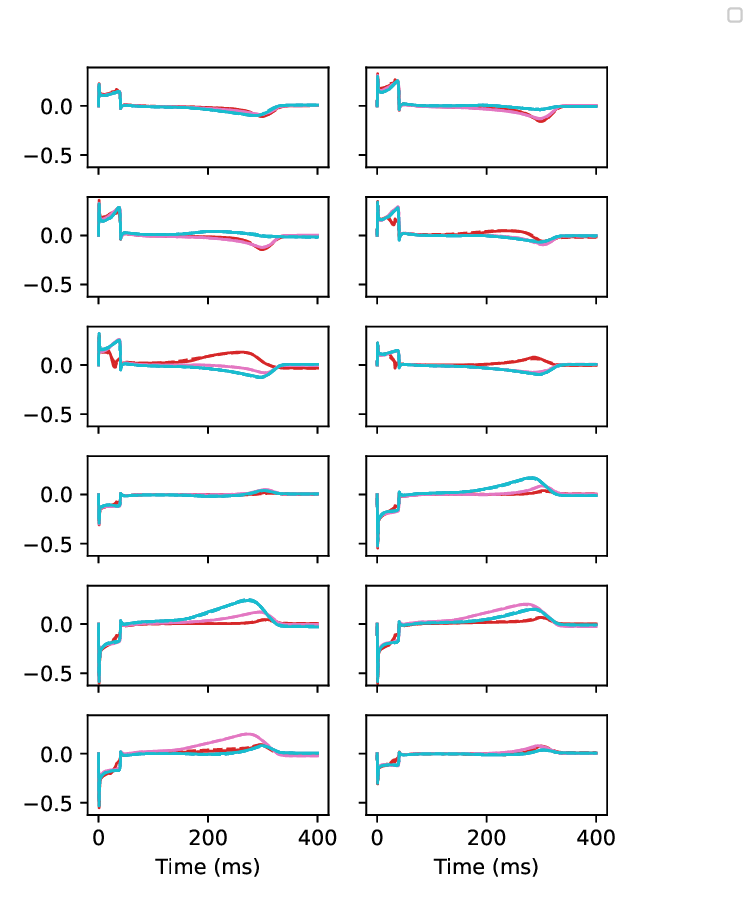}
        \caption{Forward problem. Each plot represents potential vs time for a lead for the pseudo-ECG solution.}
        \label{fig:ischemic2d_forward}
    \end{subfigure}
    \hfill
    \begin{subfigure}{0.482\textwidth}
        \centering
        
        \begin{subfigure}{\textwidth}
            \centering
            \includegraphics[width=\linewidth,
            trim={1.0cm 3.1cm 1.2cm 3.4cm},clip]
            {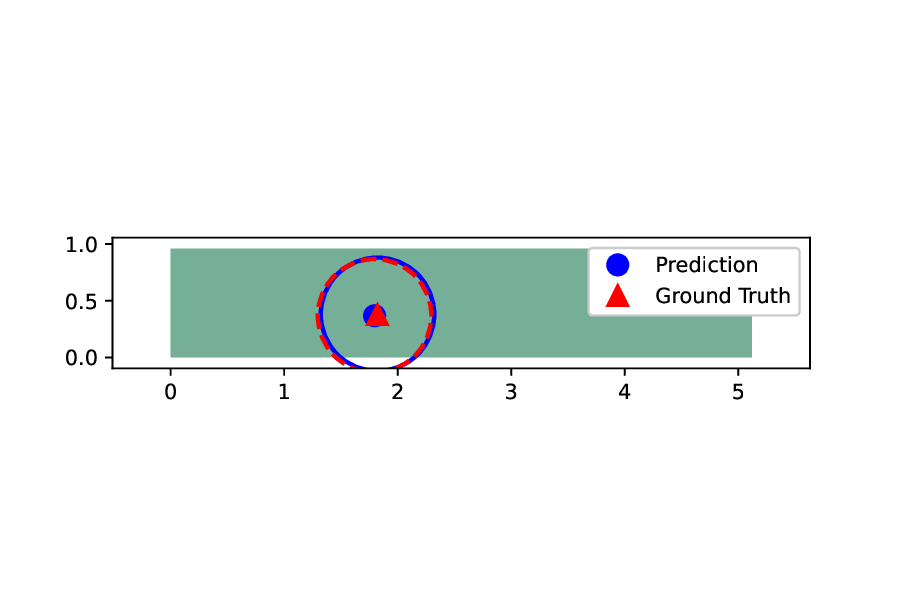}
        \end{subfigure}
        
        
        \begin{subfigure}{\textwidth}
            \centering
            \includegraphics[width=\linewidth,
            trim={1.0cm 3.1cm 1.2cm 3.4cm},clip]
            {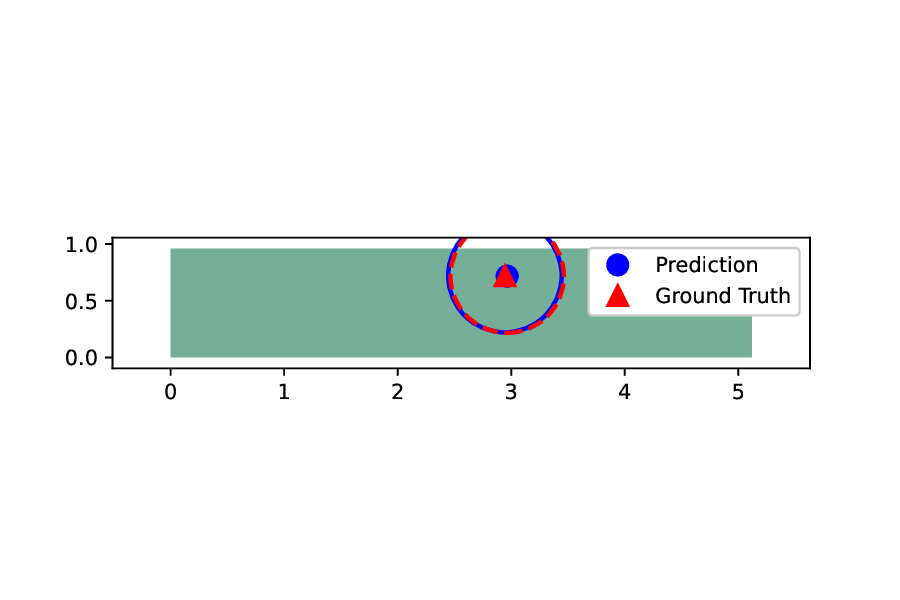}
        \end{subfigure}
        
        
        \begin{subfigure}{\textwidth}
            \centering
            \includegraphics[width=\linewidth,
            trim={1.0cm 3.1cm 1.2cm 3.4cm},clip]
            {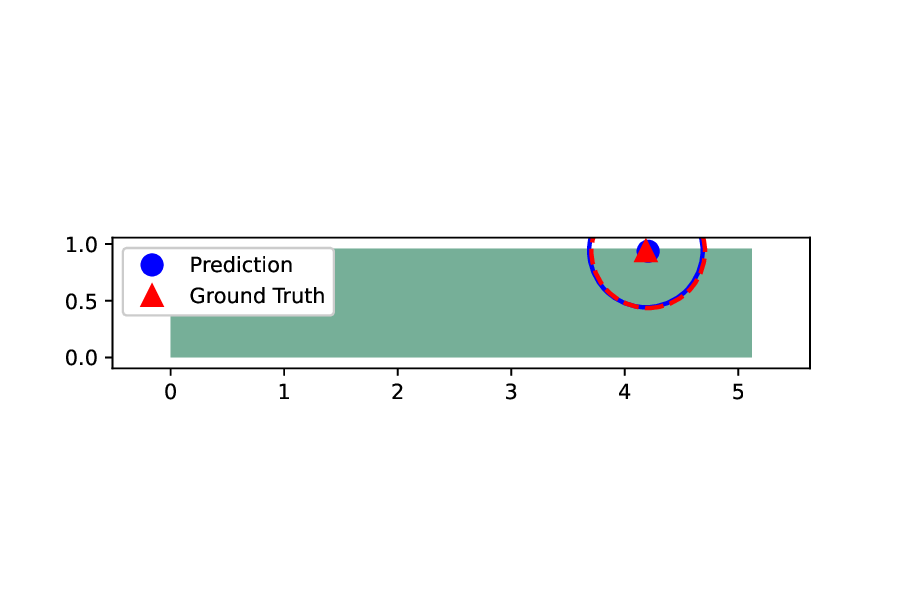}
        \end{subfigure}
        
        
        \begin{subfigure}{\textwidth}
            \centering
            \includegraphics[width=\linewidth,
            trim={1.0cm 3.1cm 1.2cm 3.4cm},clip]
            {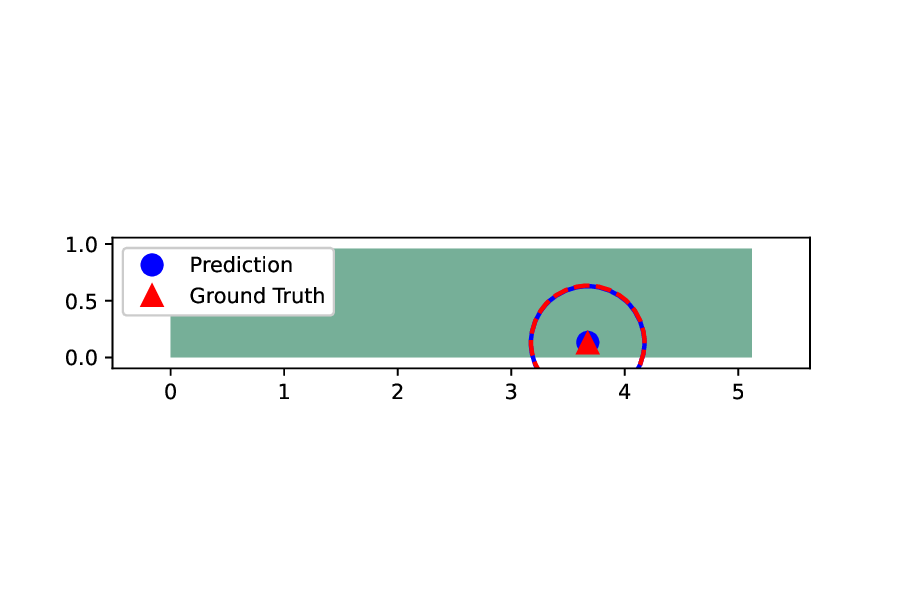}
        \end{subfigure}

        \vspace{0.53cm}
        
        \caption{Inverse problem reconstructions corresponding to the forward solutions shown on the left.}
        \label{fig:ischemic2d_inverse}
    \end{subfigure}
    
    \caption{
    2d test case for identification of the ischemic region. 
    \textbf{Left (color online):} Forward problem solutions (pseudo-ECGs generated by the ROM for different parameter instances). Different colors indicate different cases. Dashed lines indicates ground truth solution, while continuous lines indicate ML reconstructions. Results show good correspondence between ground truth and neural network approximation.
    \textbf{Right:} Corresponding inverse reconstructions obtained from those forward-generated signals. 
    Each inverse solution is computed using the pseudo-ECG displayed on the left as observed data.
    }
    \label{fig:ischemic2d_combined}
\end{figure}

In this section, the inverse problem reduces to estimating the coordinates of the centroid of the ischemic region, while the radius is assumed to be known and fixed.
Figure~\ref{fig:ischemic2d_combined} shows results of both forward and inverse problems. 
On the left, the pseudo-ECG signals generated by the surrogate model are shown for different centroid instances, proving an excellent agreement with the corresponding high-fidelity solutions. 
On the right, the inverse reconstructions highlight the capability of the proposed framework to accurately recover the location and extent of the ischemic regions from the observed pseudo-ECG data. 
The reconstructions are consistent with the forward simulations, confirming the robustness of the approach.

The architectural choices and optimization parameters used for both the forward and inverse models are reported in Table~\ref{tab:arch_params_2d_isch}. 
The forward surrogate is characterized by a moderately deep reconstruction network and a compact latent dynamics.
The inverse problem is formulated in a low-dimensional parameter space corresponding to the centroid coordinates $(x,y)$ of the ischemic region. The optimization is therefore computationally short and efficient as shown in Table \ref{tab:training_stats_2d_isch}.

Accuracy of the forward surrogate model varies together with the dimension of the latent space (cf. Table~\ref{tab:latent_states_2d_ischemic}). 
We observe that increasing the number of latent states generally improves the predictive performance, as reflected by the decrease in normalized RMSE and Pearson dissimilarity. 
In particular, the configuration with 24 latent states provides the best trade-off between accuracy and computational cost, achieving the lowest error metrics while maintaining a reasonable memory footprint and training time.

We fix the number of latent states to 24 for the solution of the inverse problem. 
Results in Table~\ref{tab:training_stats_2d_isch} indicate a low minimum loss and a satisfactory mean error over the test dataset, with moderate variability across samples. 
Moreover, the computational cost for the inverse problem remains contained, with an average optimization time below two minutes per instance and a memory usage compatible with standard CPU resources. 
Overall, these results confirm the effectiveness of the proposed data-driven strategy for the identification of ischemic regions from pseudo-ECG measurements.

\begin{table}[ht]
\centering
\begin{tblr}{
  width = \textwidth,
  colspec = {l c c},
  row{1} = {font=\bfseries},
}
\hline
Parameter & Forward model & Inverse model \\
\hline

Temporal step $dt$ & $6\times10^{-3}$ & -- \\

Dynamic network layers & $[4,\,8]$ & -- \\

Reconstruction network layers & $[17,\,17,\,17,\,17,\,17]$ & -- \\

Regularization parameter $\alpha$ & $4.7\times10^{-3}$ & -- \\

Number of subdomains & -- & $40 \times 20$ \\

Adam epochs (stage 1 / stage 2) & $300 / 300$ & $20 / 20$ \\

Adam learning rate (stage 1) & $1\times10^{-2}$ & $1\times10^{-2}$ \\

Adam learning rate (stage 2) & $1\times10^{-3}$ & $1\times10^{-3}$ \\

BFGS epochs & $15{,}000$ & $20$ \\
\hline
\end{tblr}
\caption{
Architecture and optimization parameters for the forward and inverse models for the 2d case for the identification of the ischemic region.
Here, $dt$ denotes the temporal step size of the forward Euler scheme used in the dynamics network,
the dynamic and reconstruction network layers specify the number of neurons per hidden layer,
$\alpha$ is the regularization parameter weighting the training loss,
Adam epochs and learning rates correspond to the two-stage Adam optimization procedure,
BFGS epochs indicate the number of iterations of the quasi-Newton optimizer used for fine-tuning,
while the number of subdomains defines the spatial partitioning adopted for the optimization of the inverse problem.
}
\label{tab:arch_params_2d_isch}
\end{table}

\begin{table}[ht]
\centering
\resizebox{\textwidth}{!}{  
\begin{tblr}{
  colspec = {c c c c c c c},
  row{1} = {font=\bfseries},
}
\hline
Latent states 
& Training loss 
& Validation loss 
& Normalized RMSE 
& Pearson dissimilarity 
& CPU memory [GB] 
& Training time [h] \\
\hline
8  & $3.903\times10^{-4}$ & $3.770\times10^{-4}$ & $8.919\times10^{-3}$ & $7.572\times10^{-3}$ & 5.16 & 16.87 \\
12 & $1.076\times10^{-3}$ & $8.652\times10^{-4}$ & $1.376\times10^{-2}$ & $1.816\times10^{-2}$ & 5.53 & 10.11\\
16 & $4.654\times10^{-4}$ & $3.568\times10^{-4}$ & $7.962\times10^{-3}$ & $6.031\times10^{-3}$ & 6.22 & 15.46 \\
18 & $3.361\times10^{-4}$ & $1.632\times10^{-4}$ & $7.801\times10^{-3}$ & $5.796\times10^{-3}$ & 5.75 & 15.58  \\
20 & $4.229\times10^{-4}$ & $2.334\times10^{-4}$ & $7.009\times10^{-3}$ & $4.678\times10^{-3}$ & 6.18 & 16.01 \\
24 & $3.506\times10^{-4}$ & $1.635\times10^{-4}$ & $5.816\times10^{-3}$ & $3.218\times10^{-3}$ & 6.58 & 15.99 \\
\hline
\end{tblr}
}
\caption{Performance metrics for the 2d case for the identification of the ischemic region with varying numbers of latent states (forward problem approximation).}
\label{tab:latent_states_2d_ischemic}
\end{table}

\begin{table}[ht]
\centering
\begin{tblr}{
  colspec = {l l},
  row{1} = {font=\bfseries},
}
\hline
Metric & Value \\
\hline
Minimum loss & $7.46\times10^{-8}$ \\
Mean loss (test dataset) & $3.02\times10^{-3} \pm 1.59\times10^{-2}$ \\
Maximum CPU memory usage & $10.69~\mathrm{GB}$ \\
Mean optimization time per data & $96~\mathrm{s}$ \\
\hline
\end{tblr}
\caption{Summary statistics of the inverse optimization loop for the 2d case for the identification of the ischemic region.}
\label{tab:training_stats_2d_isch}
\end{table}
\subsection{ECG-based ischemic region identification with variable radius (2d)}

\begin{table}[ht]
\centering
\begin{tblr}{
  width = \textwidth,
  colspec = {l c c},
  row{1} = {font=\bfseries},
}
\hline
Parameter & Forward model & Inverse model \\
\hline
Temporal step $dt$ & $0.12$ & -- \\

Dynamic network layers & $[6,\,6]$ & -- \\

Reconstruction network layers & $[20,\,20,\,20,\,20]$ & -- \\

Regularization parameter $\alpha$ & $4.7\times10^{-3}$ & -- \\

Number of subdomains & -- & $40 \times 20$ \\

Adam epochs (stage 1 / stage 2 / stage 3) & $200 / 200 / 200$ & $20 / 20 /$ -- \\

Adam learning rate (stage 1) & $1\times10^{-2}$ & $1\times10^{-2}$ \\

Adam learning rate (stage 2) & $1\times10^{-3}$ & $1\times10^{-3}$ \\

Adam learning rate (stage 3) & $5\times10^{-4}$ & -- \\

BFGS epochs & $10{,}000$ & $50$ \\
\hline
\end{tblr}
\caption{
Architecture and optimization parameters for the forward and inverse models for the 2d case for the identification of the ischemic region with variable radius.
Here, $dt$ denotes the temporal step size of the forward Euler scheme used in the dynamics network,
the dynamic and reconstruction network layers specify the number of neurons per hidden layer,
$\alpha$ is the regularization parameter weighting the training loss,
Adam epochs and learning rates correspond to the two-stage Adam optimization procedure,
BFGS epochs indicate the number of iterations of the quasi-Newton optimizer used for fine-tuning,
while the number of subdomains defines the spatial partitioning adopted for the optimization of the inverse problem.
}
\label{tab:arch_params_updated}
\end{table}

\begin{table}[ht]
\centering
\resizebox{\textwidth}{!}{ 
\begin{tblr}{
  colspec = {c c c c c c c},
  row{1} = {font=\bfseries},
}
\hline
Latent states 
& Training loss 
& Validation loss 
& Normalized RMSE 
& Pearson dissimilarity 
& CPU memory [GB] 
& Training time [h] \\
\hline
8  & $4.414\times10^{-4}$ & $5.470\times10^{-4}$ & $1.223\times10^{-2}$ & $1.268\times10^{-2}$ & 4.61 & 11.04 \\
12 & $4.192\times10^{-4}$ & $3.241\times10^{-4}$ & $1.028\times10^{-2}$ & $8.936\times10^{-3}$ & 5.15 & 10.84 \\
16 & $4.182\times10^{-4}$ & $2.165\times10^{-4}$ & $8.034\times10^{-3}$ & $5.449\times10^{-3}$ & 5.39 & 10.37 \\
18 & $4.566\times10^{-4}$ & $3.862\times10^{-4}$ & $1.039\times10^{-2}$ & $9.119\times10^{-3}$ & 5.40 & 10.80 \\
20 & $3.948\times10^{-4}$ & $4.217\times10^{-4}$ & $1.120\times10^{-2}$ & $1.062\times10^{-2}$ & 5.88 & 10.71 \\
24 & $4.265\times10^{-4}$ & $3.297\times10^{-4}$ & $1.005\times10^{-2}$ & $8.520\times10^{-3}$ & 6.53 & 10.86 \\
\hline
\end{tblr}
}
\caption{Performance metrics for the 2d case for the identification of the ischemic region with variable radius with varying numbers of latent states (forward problem approximation).}
\label{tab:latent_states_2d_isch_rad}
\end{table}

\begin{table}[ht]
\centering
\begin{tblr}{
  colspec = {l l},
  row{1} = {font=\bfseries},
}
\hline
Metric & Value \\
\hline
Minimum loss & $4.33\times10^{-6}$ \\
Mean loss (all data) & $1.85\times10^{-2} \pm 3.30\times10^{-3}$ \\
Maximum CPU memory usage & $12.44~\mathrm{GB}$ \\
Mean optimization time & $92~\mathrm{s}$ \\
\hline
\end{tblr}
\caption{Summary statistics of the inverse optimization loop in the 2d case for the identification of the ischemic region with variable radius.}
\label{tab:training_stats_isch}
\end{table}

Tables~\ref{tab:arch_params_updated}–\ref{tab:training_stats_isch} detail the architectural specifications, surrogate model accuracy, and inverse optimization performance for the two-dimensional variable-radius ischemic region identification problem.
In this setting, the unknown parameter vector,
\[
p = (x, y, r) \in \Omega \times \mathbb{R}^+,
\]
defines both the spatial location and the extent of a circular ischemic region characterized by modified electrophysiological properties through a reduced intracellular conductivity tensor $\mathbf{D}_i$.
Table~\ref{tab:arch_params_updated} reports the hyperparameters adopted for both the forward surrogate model and the inverse optimization procedure.
For the inverse problem, the parameter space is partitioned into $40 \times 20$ subdomains.

The accuracy of the forward surrogate as a function of the latent dimension is analyzed in Table~\ref{tab:latent_states_2d_isch_rad}.
In this case, the forward model is trained using a regularization term to the loss function presented in Equation \ref{eq:trainloss}.
This term adds the information regarding the frequency-domain.
Given predicted and target pseudo-ECG signals, the training loss reads
\[
\mathcal{L}_{\text{FFT}}(\mathcal{S}(p), \text{pECG})
= \mathcal{L}(\mathcal{S}(p), \text{pECG}) + \frac{\omega}{N_{\text{leads}}}
\sum_{i=1}^{N_{\text{leads}}}
\left\|
 \mathcal{F}(\text{pECG}_i)
-
\mathcal{F}(\widetilde{\text{pECG}}_i(\cdot;p))
\right\|^2,
\]
where $\mathcal{F}$ denotes the discrete Fourier transform applied along the temporal dimension, and $\omega=10^{-3}$ is a weighting parameter balancing the contribution of the frequency-domain term. This formulation promotes agreement between predicted and target signals not only in time but also in their spectral content, which is particularly relevant for catching the small oscillations observed in ECGs due to the varying radius of the ischemia.
Even in this case, increasing the number of latent states improves the predictive performance up to an optimal range, with the configuration employing $16$ latent states achieving the lowest validation loss, normalized RMSE, and Pearson dissimilarity. 
This indicates an improved reconstruction of both the amplitude and temporal structure of the pseudo-ECG signals. Further increases in the latent dimension do not lead to systematic accuracy gains and are accompanied by higher memory requirements, suggesting diminishing returns beyond a moderate latent dimensionality. Therefore, in this case the dimension of the latent space for surrogating the forward operator is fixed at 16.

\begin{figure}
    \centering
    \begin{subfigure}{0.46\textwidth}
        \centering
        \includegraphics[width=\linewidth,trim={0cm 0cm 1cm 0cm},clip]
        {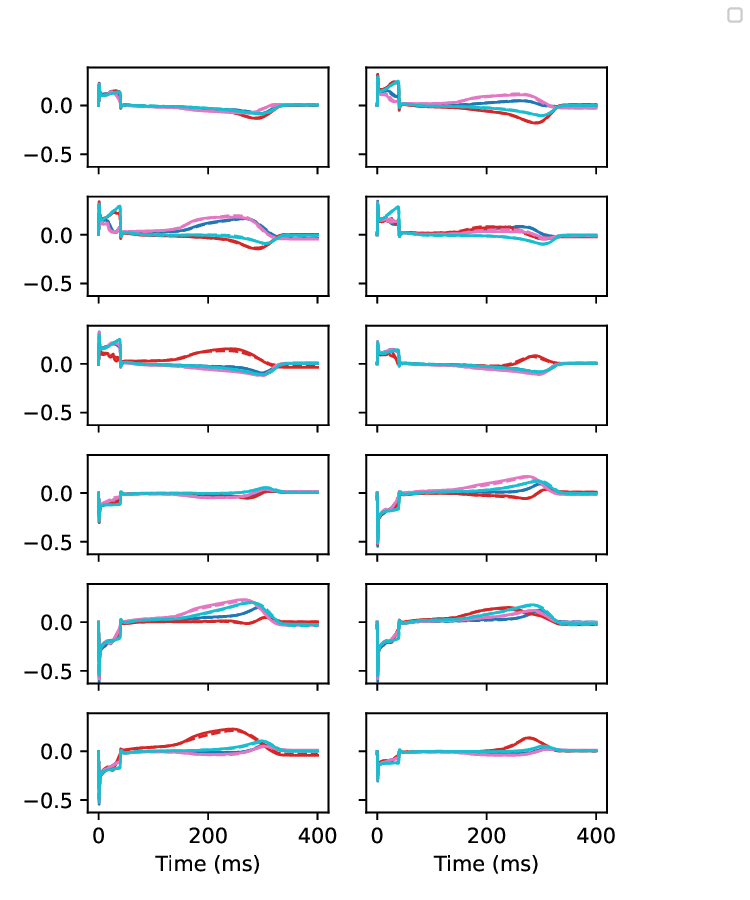}
        \caption{Forward problem. Each plot represents potential vs time for a lead for the pseudo-ECG solution.}
        \label{fig:vanilla2d_forward}
    \end{subfigure}
    \hfill
    \begin{subfigure}{0.482\textwidth}
        \centering
        
        \begin{subfigure}{\textwidth}
            \centering
            \includegraphics[width=\linewidth,
            trim={1.2cm 4.2cm 1.5cm 4.8cm},clip]
            {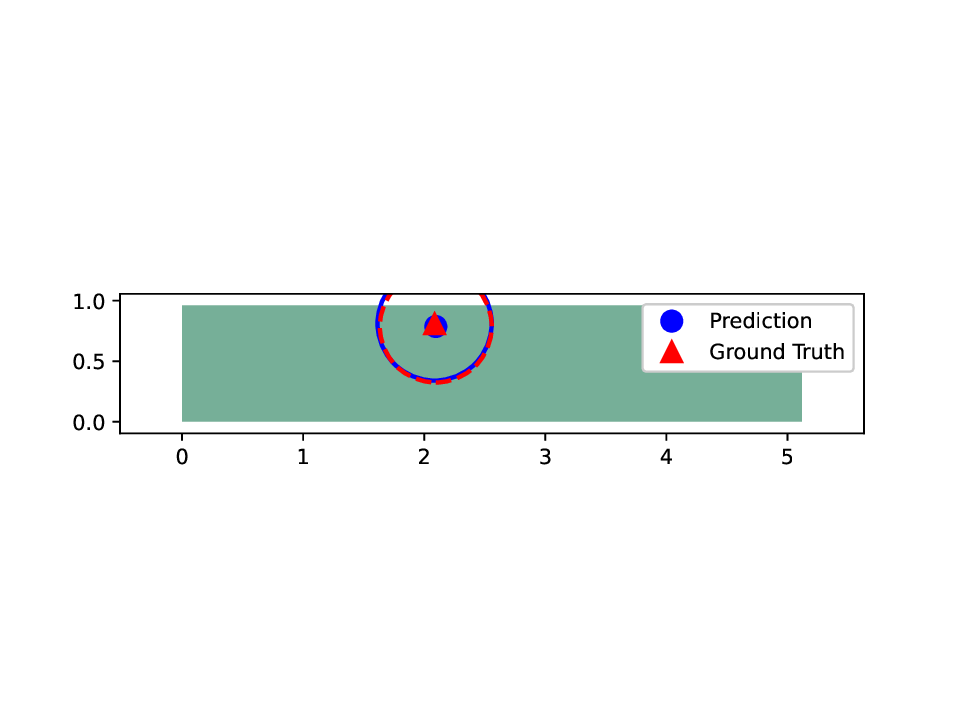}
        \end{subfigure}
        
        \vspace{0.3cm}
        
        \begin{subfigure}{\textwidth}
            \centering
            \includegraphics[width=\linewidth,
            trim={1.2cm 4.2cm 1.5cm 4.8cm},clip]
            {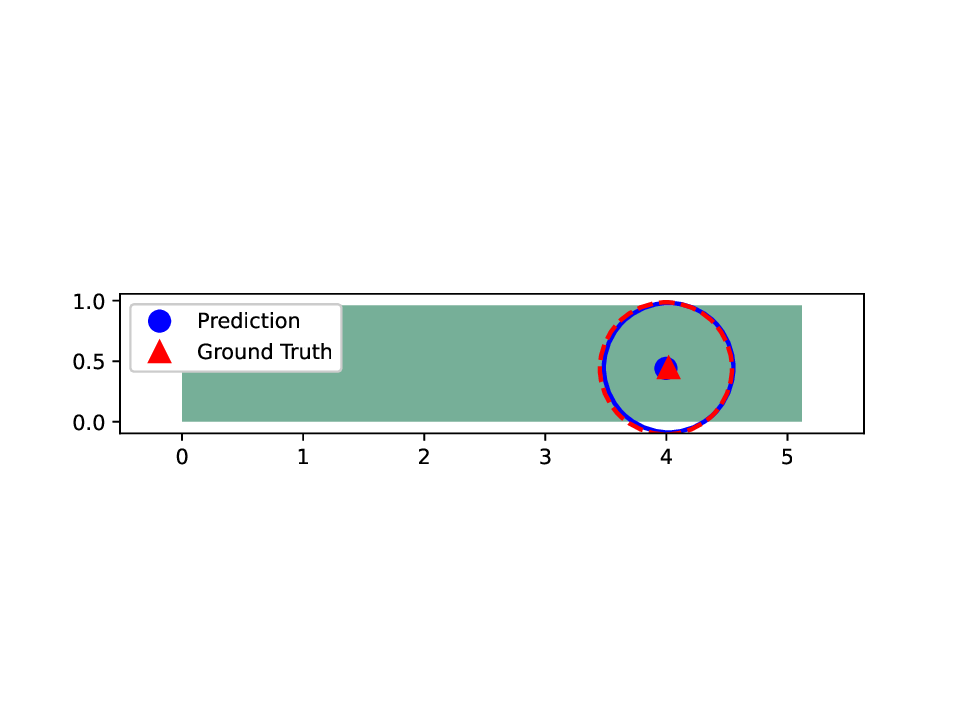}
        \end{subfigure}
        
        \vspace{0.3cm}
        
        \begin{subfigure}{\textwidth}
            \centering
            \includegraphics[width=\linewidth,
            trim={1.2cm 4.2cm 1.5cm 4.8cm},clip]
            {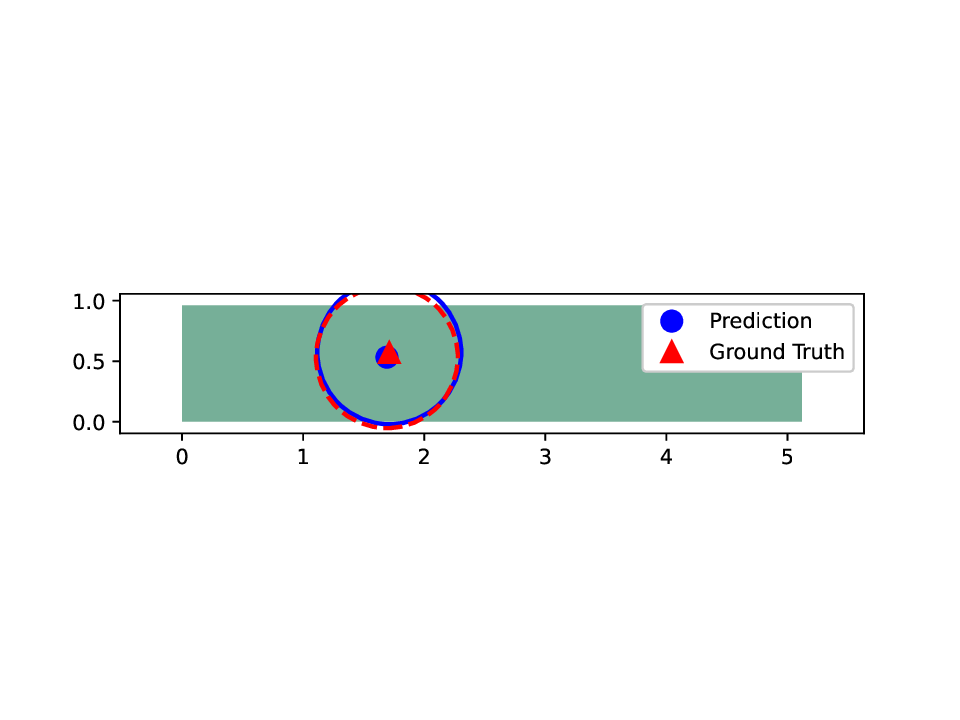}
        \end{subfigure}
        
        \vspace{0.3cm}
        
        \begin{subfigure}{\textwidth}
            \centering
            \includegraphics[width=\linewidth,
            trim={1.2cm 4.2cm 1.5cm 4.8cm},clip]
            {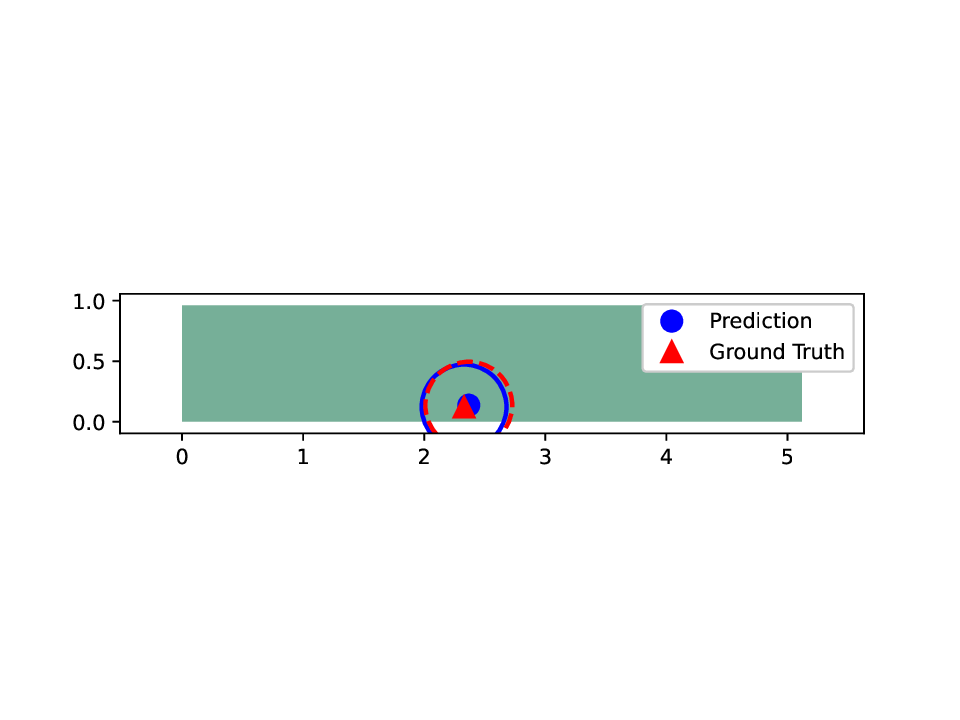}
        \end{subfigure}

        \vspace{0.65cm}
        
        \caption{Inverse problem reconstructions corresponding to the forward solutions shown on the left.}
        \label{fig:ischrad2d_inverse}
    \end{subfigure}
    
    \caption{
    2d case for the identification of the ischemic region with variable radius.
    \textbf{Left (color online):} Forward problem solutions (pseudo-ECGs generated by the ROM for different parameter instances). Different colors indicate different cases. Dashed lines indicates ground truth solution, while continuous lines indicate ML reconstructions. Results show good correspondence between ground truth and neural network approximation.
    \textbf{Right:} Corresponding inverse reconstructions obtained from those forward-generated signals. 
    Each inverse solution is computed using the pseudo-ECG displayed on the left as observed data.
    }
    \label{fig:ischrad2d_combined}
\end{figure}

Finally, Table~\ref{tab:training_stats_isch} reports aggregate statistics of the optimization loop for the inverse problem. The order of magnitude of the minimum achieved loss ($10^{-6}$) confirms the capability of the surrogate-based forward model to accurately reproduce observed pseudo-ECG signals, while the mean loss and its standard deviation indicate stable convergence across the test set.
The mean optimization time highlights the computational efficiency of the proposed framework, enabling repeated forward evaluations within a gradient-based inverse problem formulation. Overall, these results demonstrate that the surrogate model provides an accurate and efficient forward model for solving the two-dimensional ischemic region localization problem with variable radius.
Plots of the forward and the inverse problem solutions for the 2d case with ischemic region and variable radius are reported in Figure~\ref{fig:ischrad2d_combined}.

\label{sec:2d_isch_rad}
\section{Conclusions}
In this work, we propose a reduced, data-driven computational framework inspired by Latent Dynamics Networks (LDNets) \cite{reg24} for the efficient solution of inverse problems in cardiac electrophysiology.
Due to the multi-query nature of the inverse problem, the core idea is to replace the repeated evaluation of high-fidelity electrophysiological models with a surrogate model capable of accurately approximating the nonlinear mapping between low-dimensional parameters, such as initial activation sites or ischemic region descriptors, and the corresponding pseudo-ECG signals.
This approach enables a drastic reduction in computational cost while maintaining high predictive accuracy, with errors typically on the order of $10^{-4}$ in relative terms.
Furthermore, this framework can be straightforwardly extended to incorporate real clinical data, should such data become available

The end-to-end forward surrogate is designed to capture the temporal evolution of the system through a compact latent representation: in this way we are able to decouple the learning dynamics and the problem of reconstructing observable quantities.
This structure allows for efficient inference and good generalization across different parameter configurations, as confirmed by the extensive campaign of numerical experiments in both two- and three-dimensional settings.
In particular, we observed that moderately sized latent spaces provide the best trade-off between accuracy, robustness, and computational resources, avoiding overparameterization effects while retaining expressive power.
Furthermore, numerical results obtained for different benchmark problems with in-silico data, demonstrate that the proposed framework achieves accurate parameter reconstruction with moderate memory usage and reduced computational time. These features make the approach particularly appealing for other multi-query contexts (such as uncertainty quantification or optimal control) and pave the way toward real-time or near real-time applications.

The inverse problem is formulated directly in a low-dimensional parameter space, improving the stability and tractability of the optimization process.
Moreover, for the problem of identifying radius and centroid of an ischemic region we introduced a combined time–frequency loss function, which incorporates both time-domain discrepancies and spectral information via Fourier transforms. This modification proves particularly beneficial in capturing subtle dynamical features of the pseudo-ECG signals and improves the identifiability of the underlying parameters.

Convergence of the inverse solver is ensured by combining a tailored, domain-partitioned multi-start initialization with a hybrid optimization routine that leverages both Adam and second order quasi-Newton (BFGS) methods.
In more challenging scenarios, such  the identification of the stimulus in three dimensional domains, we prove that a multi-start strategy based on multiple initial guesses provides a robust compromise between computational efficiency and reconstruction accuracy.

Overall, the proposed methodology represents a significant step toward practical, noninvasive, and data-driven solutions to inverse problems in electrocardiology.
Future developments will focus on extending the framework to more realistic anatomical geometries, incorporating patient-specific data, and further improving robustness with respect to noise and model uncertainties, with the ultimate goal of supporting clinical decision-making processes such as ablation planning and the localization of arrhythmic sources.

\section*{Acknowledgements}
{\small
EC, LP and SS acknowledge the CINECA award under the ISCRA initiative, for the availability of high-performance computing resources and support (projects DDO2CARD and DDML2Car). EC, LP, SS and GZ are members of INdAM-GNCS. EC and LP have been supported by MUR (PRIN 202232A8AN\_002 and PRIN P2022B38NR\_001) funded by European Union - Next Generation EU. SS and GZ have been supported by MUR (PRIN 202232A8AN\_003 and PRIN P2022B38NR\_002) funded by European Union - Next Generation EU. GZ has received support from the project FIS, MUR, Italy 2025-2028, Project code: FIS-2023-02228, CUP: D53C24005440001, ``SYNERGIZE: Synergizing Numerical Methods and Machine Learning for a new generation of computational models''. The funders had no role in study design, data collection and analysis, decision to publish, or preparation of the manuscript.
}


\begin{thebibliography}{99} 

\bibitem{afr22} Africa, P. C. (2022). lifex: A flexible, high performance library for the numerical solution of complex finite element problems. \textit{SoftwareX, 20}, 101252.

\bibitem{asp25}
Aspri, A., Beretta, E., Francini, E., Pierotti, D., \& Vessella, S. (2025). On an inverse problem with applications in cardiac electrophysiology. \textit{Nonlinearity, 38}(4), 045014.

\bibitem{azi24} Azizzadenesheli, K., Kovachki, N., Li, Z., Liu-Schiaffini, M., Kossaifi, J., \& Anandkumar, A. (2024). Neural operators for accelerating scientific simulations and design. \textit{Nature Reviews Physics, 6}(5), 320-328.

\bibitem{bab04} Bakushinsky, A. B., \& Kokurin, M. Y. (2004). \textit{Iterative Methods for Approximate Solution of Inverse Problems}. Springer.

\bibitem{bal23} Balay, S., Abhyankar, S., Adams, M., Benson, S., Brown, J., Brune, P., ... \& Zhang, J. (2021). \textit{PETSc/TAO users manual} (ANL-21/39-Revision 3.17). Argonne National Laboratory.

\bibitem{bat24} Batlle, P., Darcy, M., Hosseini, B., \& Owhadi, H. (2024). Kernel methods are competitive for operator learning. \textit{Journal of Computational Physics, 496}, 112549.

\bibitem{ben15} Benner, P., Gugercin, S., \& Willcox, K. (2015). A survey of projection-based model reduction methods for parametric dynamical systems. \textit{SIAM Review, 57}(4), 483-531.

\bibitem{ber17} Beretta, E., Cavaterra, C., Cerutti, M. C., Manzoni, A., \& Ratti, L. (2017). An inverse problem for a semilinear parabolic equation arising from cardiac electrophysiology. Inverse Problems, 33(10), 105008.

\bibitem{ber20} Beretta, E., Cavaterra, C., \& Ratti, L. (2020). On the determination of ischemic regions in the monodomain model of cardiac electrophysiology from boundary measurements. Nonlinearity, 33(11), 5659-5685.

\bibitem{bru16} Brunton, S. L., Proctor, J. L., \& Kutz, J. N. (2016). Discovering governing equations from data by sparse identification of nonlinear dynamical systems. \textit{Proceedings of the National Academy of Sciences, 113}(15), 3932-3937.

\bibitem{cen24} Centofanti, E., \& Scacchi, S. (2024). A comparison of algebraic multigrid bidomain solvers on hybrid CPU–GPU architectures. \textit{Computer Methods in Applied Mechanics and Engineering, 423}, 116875.

\bibitem{col14} Colli Franzone, P., Pavarino, L. F., \& Scacchi, S. (2014). \textit{Mathematical cardiac electrophysiology} (Vol. 13). Springer.

\bibitem{dok26} Dokuchaev, A., Bonizzoni, F., Pagani, S., Regazzoni, F., \& Pezzuto, S. (2026). Learning geometry-dependent lead-field operators for forward ECG modeling. arXiv preprint arXiv:2602.22367.

\bibitem{dut17} Dutta, S., Mincholé, A., Quinn, T. A., \& Rodriguez, B. (2017). Electrophysiological properties of computational human ventricular cell action potential models under acute ischemic conditions. \textit{Progress in Biophysics and Molecular Biology, 129}, 40-52.

\bibitem{fal06} Falgout, R. D., Jones, J. E., \& Yang, U. M. (2006). The design and implementation of hypre, a library of parallel high performance preconditioners. In A. M. Bruaset \& A. Tveito (Eds.), \textit{Numerical solution of partial differential equations on parallel computers} (pp. 267–294). Springer.

\bibitem{ges83} Geselowitz, D. B., \& Miller III, W. T. (1983). A bidomain model for anisotropic cardiac muscle. \textit{Annals of Biomedical Engineering, 11}(3), 191-206.

\bibitem{gos23} Goswami, S., Bora, A., Yu, Y., \& Karniadakis, G. E. (2023). Physics-informed deep neural operator networks. In \textit{Machine learning in modeling and simulation: Methods and applications} (pp. 219–254). Springer International Publishing.

\bibitem{gra22} Grandits, T., Pezzuto, S., \& Plank, G. (2022). Smoothness and continuity of cost functionals for ECG mismatch computation. \textit{IFAC-PapersOnLine, 55}(20), 181-186.

\bibitem{gra25} Grandits, T., Gillette, K., Plank, G., \& Pezzuto, S. (2025). Accurate and efficient cardiac digital twin from surface ECGs: Insights into identifiability of ventricular conduction system. \textit{Medical Image Analysis}.

\bibitem{izh06} Izhikevich, E. M., \& FitzHugh, R. (2006). FitzHugh-Nagumo model. \textit{Scholarpedia, 1}(9), 1349.

\bibitem{li20} Li, Z., Kovachki, N., Azizzadenesheli, K., Liu, B., Bhattacharya, K., Stuart, A., \& Anandkumar, A. (2020). Fourier neural operator for parametric partial differential equations. \textit{arXiv preprint arXiv:2010.08895}.

\bibitem{lin22} Linot, A. J., \& Graham, M. D. (2022). Data-driven reduced-order modeling of spatiotemporal chaos with neural ordinary differential equations. \textit{Chaos: An Interdisciplinary Journal of Nonlinear Science, 32}(7).

\bibitem{liu23} Liu, B., Ocegueda, E., Trautner, M., Stuart, A. M., \& Bhattacharya, K. (2023). Learning macroscopic internal variables and history dependence from microscopic models. \textit{Journal of the Mechanics and Physics of Solids, 178}, 105329.

\bibitem{lu21} Lu, L., Jin, P., Pang, G., Zhang, Z., \& Karniadakis, G. E. (2021). Learning nonlinear operators via DeepONet based on the universal approximation theorem of operators. \textit{Nature Machine Intelligence, 3}(3), 218-229.

\bibitem{mau21} Maulik, R., Lusch, B., \& Balaprakash, P. (2021). Reduced-order modeling of advection-dominated systems with recurrent neural networks and convolutional autoencoders. \textit{Physics of Fluids, 33}(3).

\bibitem{nie07} Nielsen, B. F., Lysaker, M., \& Tveito, A. (2007). On the use of the resting potential and level set methods for identifying ischemic heart disease: An inverse problem. Journal of Computational Physics, 220(2), 772-790.

\bibitem{nie09} Nielsen, B. F., Cai, X., Sundnes, J., \& Tveito, A. (2009). Towards a computational method for imaging the extracellular potassium concentration during regional ischemia. Mathematical biosciences, 220(2), 118-130.

\bibitem{noc06} Nocedal, J., \& Wright, S. J. (2006). \textit{Numerical optimization}. Springer.

\bibitem{pal14} Palamara, S., Vergara, C., Catanzariti, D., Faggiano, E., Pangrazzi, C., Centonze, M., ... \& Quarteroni, A. (2014). Computational generation of the Purkinje network driven by clinical measurements: The case of pathological propagations. \textit{International Journal for Numerical Methods in Biomedical Engineering, 30}(12), 1558-1577.

\bibitem{pul10} Pullan, A. J., Cheng, L. K., Nash, M. P., Ghodrati, A., MacLeod, R., \& Brooks, D. H. (2010). The inverse problem of electrocardiography. In \textit{Comprehensive electrocardiology} (pp. 299–344). Springer.

\bibitem{raj20} Rajendra, P., \& Brahmajirao, V. (2020). Modeling of dynamical systems through deep learning. \textit{Biophysical Reviews, 12}(6), 1311-1320.

\bibitem{reg19} Regazzoni, F., Dede, L., \& Quarteroni, A. (2019). Machine learning for fast and reliable solution of time-dependent differential equations. \textit{Journal of Computational Physics, 397}, 108852.

\bibitem{reg24} Regazzoni, F., Pagani, S., Salvador, M., Dede’, L., \& Quarteroni, A. (2024). Learning the intrinsic dynamics of spatio-temporal processes through Latent Dynamics Networks. \textit{Nature Communications, 15}(1), 1834.

\bibitem{ruu08} Ruud, T. S., Nielsen, B. F., Lysaker, M., \& Sundnes, J. (2008). A computationally efficient method for determining the size and location of myocardial ischemia. IEEE Transactions on Biomedical Engineering, 56(2), 263-272.

\bibitem{sca23} Scacchi, S., Colli Franzone, P., Pavarino, L. F., Gionti, V., \& Storti, C. (2023). Epicardial dispersion of repolarization promotes the onset of reentry in Brugada syndrome: A numerical simulation study. \textit{Bulletin of Mathematical Biology, 85}(3), 22.

\bibitem{sit20} Sitzmann, V., Martel, J., Bergman, A., Lindell, D., \& Wetzstein, G. (2020). Implicit neural representations with periodic activation functions. \textit{Advances in Neural Information Processing Systems, 33}, 7462-7473.

\bibitem{sun06} Sundnes, J., Lines, G. T., Cai, X., Nielsen, B. F., Mardal, K.-A., \& Tveito, A. (2007). \textit{Computing the electrical activity in the heart} (Vol. 1). Springer.

\bibitem{ten06} Ten Tusscher, K. H., \& Panfilov, A. V. (2006). Cell model for efficient simulation of wave propagation in human ventricular tissue under normal and pathological conditions. \textit{Physics in Medicine \& Biology, 51}(23), 6141.

\bibitem{tur24} Turisini, M., Cestari, M., \& Amati, G. (2024). LEONARDO: A pan-European pre-exascale supercomputer for HPC and AI applications. \textit{Journal of Large-Scale Research Facilities, 9}(1), 1–16.

\bibitem{ver16} Vergara, C., Lange, M., Palamara, S., Lassila, T., Frangi, A. F., \& Quarteroni, A. (2016). A coupled 3D–1D numerical monodomain solver for cardiac electrical activation in the myocardium with detailed Purkinje network. \textit{Journal of Computational Physics, 308}, 218–238.

\bibitem{wan14} Wang, W., Huang, Y., Wang, Y., \& Wang, L. (2014). Generalized autoencoder: A neural network framework for dimensionality reduction. In \textit{Proceedings of the IEEE Conference on Computer Vision and Pattern Recognition Workshops} (pp. 490-497).

\bibitem{zap25} Zappon, E., Azzolin, L., Gsell, M. A., Thaler, F., Prassl, A. J., Arnold, R., ... \& Plank, G. (2025). An efficient end-to-end computational framework for the generation of ECG calibrated volumetric models of human atrial electrophysiology. Medical image analysis, 103822.

\bibitem{zia25} Ziarelli, G., Pagani, S., Parolini, N., Regazzoni, F., \& Verani, M. (2025). A model learning framework for inferring the dynamics of transmission rate depending on exogenous variables for epidemic forecasts. \textit{Computer Methods in Applied Mechanics and Engineering, 437}, 117796.

\bibitem{zia26} Ziarelli, G., Centofanti, E., Parolini, N., Scacchi, S., Verani, M., \& Pavarino, L. F. (2026). Learning cardiac activation and repolarization times with operator learning. \textit{PLOS Computational Biology, 22}(1), e1013920.
\end{thebibliography}
\end{document}